\documentclass[a4paper,11pt]{article}
\usepackage{bm,mathrsfs,amsmath,amsthm,amssymb,eepic,amscd,longtable,array,graphicx}
\newcommand{\nc}{\newcommand}
\nc{\rnc}{\renewcommand}
\nc{\nn}{\nonumber}
\nc{\der}{{\partial}}
\rnc{\Im}{{\rm{Im}\,}}
\rnc{\Re}{{\rm{Re}\,}}
\nc{\db}{\displaybreak[0]\\}
\nc{\bra}{\langle}
\nc{\ket}{\rangle}
\nc{\bs}{\boldsymbol}

\newtheorem{theorem}{Theorem}[section]
\newtheorem{lemma}[theorem]{Lemma}

\newtheorem{proposition}[theorem]{Proposition}
\newtheorem{corollary}[theorem]{Corollary}

\theoremstyle{definition}
\newtheorem{definition}[theorem]{Definition}

\numberwithin{equation}{section}

\numberwithin{equation}{section}

\begin{document}%
%
\title{Dual wavefunction of the symplectic ice}

\author{
Kohei Motegi \thanks{E-mail: kmoteg0@kaiyodai.ac.jp}
\\\\
{\it Faculty of Marine Technology, Tokyo University of Marine Science and Technology,}\\
 {\it Etchujima 2-1-6, Koto-Ku, Tokyo, 135-8533, Japan} \\
\\\\
\\
}

\date{\today}

\maketitle

\begin{abstract}
The wavefunction of the free-fermion six-vertex model
was found to give a natural realization of the Tokuyama combinatorial formula
for the Schur polynomials by Bump-Brubaker-Friedberg.
Recently, we studied the correspondence between the
dual version of the wavefunction and the Schur polynomials,
which gave rise to another combinatorial formula.
In this paper, we extend the analysis to the reflecting boundary condition,
and show the exact correspondence between
the dual wavefunction and the symplectic Schur functions.
This gives a dual version of the integrable model realization of
the symplectic Schur functions by Ivanov.
We also generalize to the correspondence between the
wavefunction, the dual wavefunction of the six-vertex model
and the factorial symplectic Schur functions
by the inhomogeneous generalization of the model.
\end{abstract}
{\bf Mathematics Subject Classification.}
05E05, 05E10, 16T25, 16T30, 17B37. \\
\\
{\bf Keywords.} Integrable lattice models, Yang-Baxter equation,
Symmetric functions, Combinatorial representation theeory.

\section{Introduction}
The study of the connections between
integrable lattice models \cite{Bethe,FST,Baxter,KBI}
and combinatorial representation theory of symmetric polynomials
is an active area of research in these years.

The main actor is the wavefunction in this research,
which is constructed from the $R$-matrices satisfying the Yang-Baxter relation.
The wavefunction of
the most famous six-vertex model \cite{Lieb,Reshetikhin}
whose $R$-matrix comes from the Drinfeld-Jimbo 
quantum group \cite{Dr,J} $U_q(sl_2)$ representation,
and its $q=0$ five-vertex model degeneration,
has turned out to give a representation
of the Grothendieck polynomials and its quantum group deformation
(see \cite{MS,MS2,Borodin,BP1,WZ,BW,BWZ,Korff,GK2} for example).
Based on this correspondence, various algebraic combinatorial identities
have been discovered and proved. It seems that many
of the identities themselves are hard
to be discovered without the power of quantum integrability.
These developments may shed new light in the world of
Schubert calculus \cite{LS,FK,Buch,IN,IS,Mc,KirillovSigma,Matsumura}.
For example, the notion of excited Young diagram \cite{IN}
in the field of Schubert calculus is essentially
equivalent to the wavefunction of certain integrable five-vertex models.
Translating into the language of quantum integrable models
can give new insights.
For example,
there is a recent development on the investigation of
the Littlewood-Richardson coefficients
from the point of view of quantum integrability \cite{WZ}.

The object treated in this paper
is based on another way of direction of the developments
uncovered by number theorists.
Bump, Brubaker and Friedberg found \cite{BBF} that
a ceratin kind of free-fermion six-vertex model
gives a natural realization of the Tokuyama formula \cite{To},
which gives a deformation of the Weyl character formula.
The free-fermion six-vertex model can be regarded as
a gauge-deformed version of the Felderhof free-fermion model \cite{Felderhof},
whose underlying quantum group symmetry
can be explained either as an exotic roots of unity
finite-dimensional highest weight representation
\cite{Mu,DA}.
Another explanation is that the free-fermion
six-vertex model corresponds to the simplest case
of the Perk-Schultz model \cite{PS} which is a representation
of quantum supergroup.
The latter formulation recently gave rise to the
Yang-Baxter equation \cite{BBB} for the metaplectic ice \cite{Meta}.
As for the domain wall boundary partition function
which is a special class of partition functions,
it was evaluated in \cite{ZZ,FCWZ}
by using the Izergin-Korepin technique in the past
and factorization phenomena was observed.

However, it was only found in recent years that
the free-fermion model has rich mathematical structures
related with the combinatorial representation theory
of Schur polynomials.
One of the striking facts found \cite{BBF} was that
the Tokuyama formula \cite{To}, which is a
one-parameter deformation of the 
Weyl character formula,
is naturally realized as the wavefunction of
the free-fermion model.
The wavefunction is the most fundamental object
in the statistical physical aspects of quantum integrable models,
since it becomes the Bethe eigenvectors of the
corresponding one-dimensional spin chain
when the Bethe ansatz equation is imposed
on the spectral parameters.

The Tokuyama formula for the Schur polynomials can be understood
as a consequence of the evaluation of the wavefunction
in two ways.
One by expressing it as a product of a one-parameter deformation
of the Vandermonde determinant and the Schur polynomials,
and another one by making a microscopic analysis and
deirve an expression using the strict Gelfand-Tsetlin pattern.
The Tokuyama formula is a consequence of the two evaluations
for the same object.
This understanding \cite{BBF} opened a new doorway
to the combinatorial representation theory of symmetric polynomials
via the free-fermion model.
One of the progresses after this work is the construction of
variations of the Tokuyamya-type formulas
by changing boundary conditions.
Several Tokuyama-type formulas for types other than type $A$
was obtained by Okada \cite{OkTo} and Hamel-King \cite{HK1,HK2} earlier,
and recently by using methods of analytic number theory \cite{FZ}
and non-intersecting lattice paths \cite{HKnew,HK,HKFPSAC}.

They were investigated from the viewpoint of
quantum integrability
in \cite{Iv,BBCG,Tabony,BS}
where local objects and relations such as the $L$-operators,
$K$-matrix and the Yang-Baxter relation are extensively used.
For example, the correspondence between the
wavefunction under the half-turn boundary condition
and the symplectic Schur functions was obtained by Ivanov in \cite{Iv}.

On the other hand, we studied the dual wavefunction
of the free-fermion model in a recent paper \cite{LMP}.
We gave the exact correspondence between the
dual wavefunction and the Schur polynomials.
which includes the special case $t=1$ of the deformation parameter
\cite{BBF,BMN}, which was obtained by transforming the original wavefuncion
to the dual wavefunction by symmetry arguments.
We gave two proofs for the correspondence.
One proof used transformation of the statement of the theorem
to an equivalent statement so that
one can use the arguments given in \cite{BBF}.
Another proof is a modern statistical mechanical approach,
which combines the matrix product method \cite{GMmat,KM}
and the Izergin-Korepin method of analysis
on the domain wall boundary partition function \cite{Ko,Iz}.
The exact correspondence with the Schur polynomials,
together with a microscopic analysis of the dual wavefunction
gave rise to a dual version of the Tokuyama-type formula
for the Schur polynomials.

In this paper, we combine these two directions of progresses,
and extend the study of the dual wavefunction to the free-fermion model
under the reflecting boundary condition.
We give the exact correspondence between the
dual wavefunction and the symplectic Schur functions.
We prove the correspondence by extending the argument used in
the first proof in \cite{LMP} to the reflecting boundary condition,
so that one can use the arguments by Brubaker-Bump-Friedberg \cite{BBF}
and Ivanov \cite{Iv}.
This gives a dual version of type $C$ Tokuyama formula.
We also generalize the Theorems by Ivanov and the main result
in this paper to give the exact correspondence between
the wavefunction, dual wavefunction and the
factorial symplectic Schur functions
by using two types of $L$-operators and the $K$-matrix of the inhomogeneous
six-vertex model.
This is a symplectic analogue of the work of Bump-McNamara-Nakasuji
\cite{BMN} which they established the Tokuyama-type formulas
for factorial Schur functions.
Recently, the Tokuyama-type formulas for factorial characters
were obtained by using non-intersecting lattice paths by
Hamel-King \cite{HK,HKFPSAC} which is worthwhile to
investigate the relation and the results in this paper in the future.

This paper is organized as follows.
We introduce the free-fermion model in section 2 and
review the relation between the wavefunction
under the reflecting boundary condition and
the symplectic Schur functions in section 3.
In section 4, we introduce the dual wavefunction,
and prove the relation with the symplectic Schur functions,
which can be regarded as a dual version of type $C$ Tokuyama formula.
In section 5, we generalize the correspondence and give the
exact relation between the wavefuncion, dual wavefunction
of the inhomogeneous free-fermion model
and the factorial symplectic Schur functions.
Section 6 is devoted to the conclusion.

\section{Free-fermion model}
We introduce the free-fermion model in this section,
and review
the results on the relation between the wavefunction and the
symplectic Schur functions in the next section.
We use the $L$-operator in \cite{BBF} which is
best suited for the study of the combinatorics of the Schur polynomials,
since the Tokuyama formula is exactly realized as the wavefunction
constructed from this $L$-operator.
More generic or gauge-transformed
ones can be found in \cite{Mu,DA,FCWZ} for example.
We also use the terminology of the quantum inverse scattering method
or the algebraic Bethe ansatz,
which is one of the most fundamental methods for the analysis
of quantum integrable models.

The most fundamental objects in integrable lattice models
are the $R$-matrix and the $L$-operator.
For the case of the free-fermion model we consider,
the $R$-matrix is given by
\begin{eqnarray}
R_{ab}(z,t)=\left( 
\begin{array}{cccc}
1+tz^{-1} & 0 & 0 & 0 \\
0 & t(1-z^{-1}) & t+1 & 0 \\
0 & (t+1)z^{-1} & z^{-1}-1 & 0 \\
0 & 0 & 0 & z^{-1}+t
\end{array}
\right), \label{rmatrix}
\end{eqnarray}
acting on the tensor product $W_a \otimes W_b$
of the complex two-dimensional space $W_a$.
Let us denote the orthonormal basis of $W_a$ and its dual as
$\{|0 \rangle_a, |1 \rangle_a \}$ and $\{{}_a \langle 0|, {}_a \langle 1|\}$,
and the matrix elements of the $R$-matrix as
$
{}_a \langle \gamma | {}_b \langle \delta | R_{a b}(z,t)
|\alpha \rangle_a | \beta \rangle_b=[R_{a b}(z,t)]_{\alpha \beta}^{\gamma \delta}
$. The matrix elements of the $R$-matrix are explicitly given as
\begin{align}
{}_a \langle 0| {}_b \langle 0 | R_{a b}(z,t)
|0 \rangle_a | 0 \rangle_b&=1+tz^{-1}, \\
{}_a \langle 0| {}_b \langle 1 | R_{a b}(z,t)
|0 \rangle_a | 1 \rangle_b&=t(1-z^{-1}), \\
{}_a \langle 0| {}_b \langle 1 | R_{a b}(z,t)
|1 \rangle_a | 0 \rangle_b&=t+1, \\
{}_a \langle 1| {}_b \langle 0 | R_{a b}(z,t)
|0 \rangle_a |1 \rangle_b&=(t+1)z^{-1}, \\
{}_a \langle 1| {}_b \langle 0 | R_{a b}(z,t)
|1 \rangle_a | 0 \rangle_b&=z^{-1}-1, \\
{}_a \langle 1| {}_b \langle 1 | R_{a b}(z,t)
|1 \rangle_a | 1 \rangle_b&=z^{-1}+t.
\end{align}

\begin{figure}[ht]
\includegraphics[width=12cm]{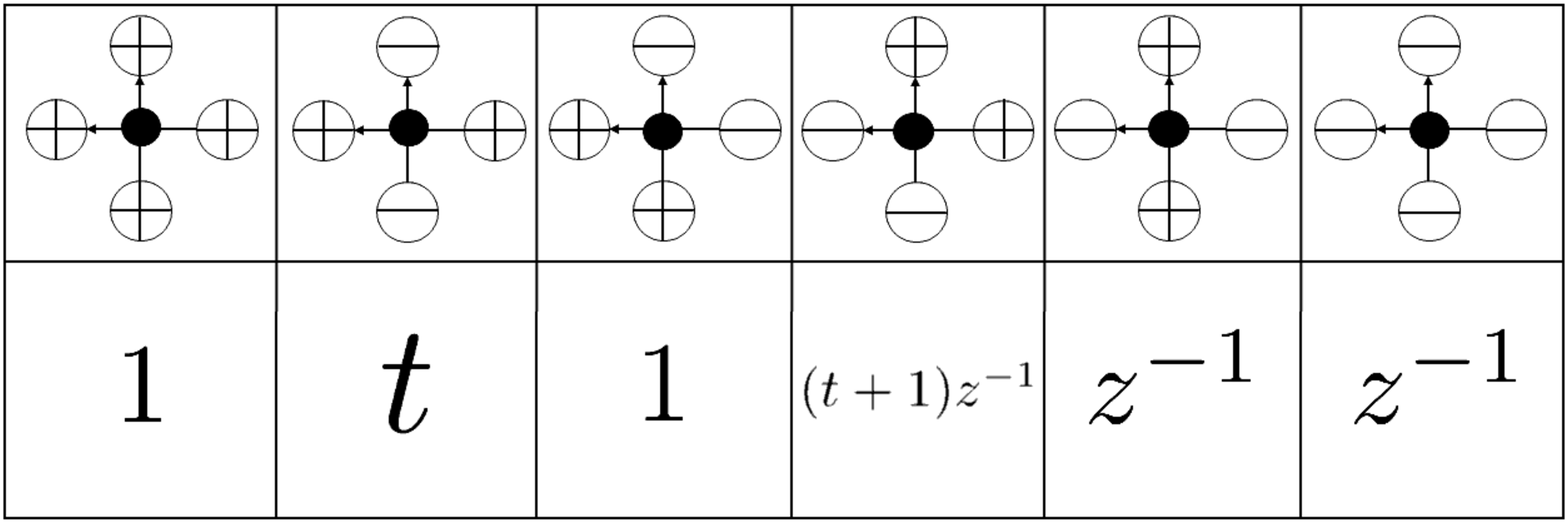}
\caption{The first $L$-operator \eqref{loperator}. The (dual) state $| 0 \rangle$ ($\langle 0 |$) is represented as $\oplus$,
while the (dual) state $| 1 \rangle$ ($\langle 1 |$)
is represented as $\ominus$, following the pictorial description of \cite{BBF}.}
\label{pictureloperator}
\end{figure}

The $L$-operator of the free-fermion model is given by
\begin{eqnarray}
L_{aj}(z,t)=\left( 
\begin{array}{cccc}
1 & 0 & 0 & 0 \\
0 & t & 1 & 0 \\
0 & (t+1)z^{-1} & z^{-1} & 0 \\
0 & 0 & 0 & z^{-1}
\end{array}
\right), \label{loperator}
\end{eqnarray}
acting on the tensor product $W_a \otimes \mathcal{F}_j$ of the space
$W_a$ and the two-dimensional Fock space at the
$j$th site $\mathcal{F}_j$.
We also denote the orthonormal basis of $\mathcal{F}_j$ and its dual as
$\{|0 \rangle_j, |1 \rangle_j \}$ and $\{{}_j \langle 0|, {}_j \langle 1|\}$,
and the matrix elements of the $L$-operator as
$
{}_a \langle \gamma| {}_j \langle \delta | L_{a j}(z,t)
|\alpha \rangle_a | \beta \rangle_j=[L_{aj}(z,t)]_{\alpha \beta}^{\gamma \delta}
$.
The matrix elements of the $L$-operator are explicitly written as
(see Figure \ref{pictureloperator} for a pictorial description)
\begin{align}
{}_a \langle 0| {}_j \langle 0 | L_{a j}(z,t)
|0 \rangle_a | 0 \rangle_j&=1, \\
{}_a \langle 0| {}_j \langle 1 | L_{a j}(z,t)
|0 \rangle_a | 1 \rangle_j&=t, \\
{}_a \langle 0| {}_j \langle 1 | L_{a j}(z,t)
|1 \rangle_a | 0 \rangle_j&=1, \\
{}_a \langle 1| {}_j \langle 0 | L_{a j}(z,t)
|0 \rangle_a |1 \rangle_j&=(t+1)z^{-1}, \\
{}_a \langle 1| {}_j \langle 0 | L_{a j}(z,t)
|1 \rangle_a | 0 \rangle_j&=z^{-1}, \\
{}_a \langle 1| {}_j \langle 1 | L_{a j}(z,t)
|1 \rangle_a | 1 \rangle_j&=z^{-1}.
\end{align}
The $R$-matrices and
the $L$-operators have origins in statistical physics,
and $| 0 \rangle$ or its dual $\langle 0|$
can be regarded as a hole state,
while $| 1 \rangle$ or its dual $\langle 1|$
can be interpretted as a particle state
from the point of view of statistical physics.
We use the terms hole states and particle states
to describe states constructed from
$| 0 \rangle$, $\langle 0|$, $| 1 \rangle$ and $\langle 1|$
from now on since they are convenient for the description
of the states.
We also remark that in the language of the quantum inverse scattering method,
the Fock spaces $W_a$ and $\mathcal{F}_j$
are usually called as the auxiliary and quantum spaces, respectively.

\begin{figure}[ht]
\includegraphics[width=12cm]{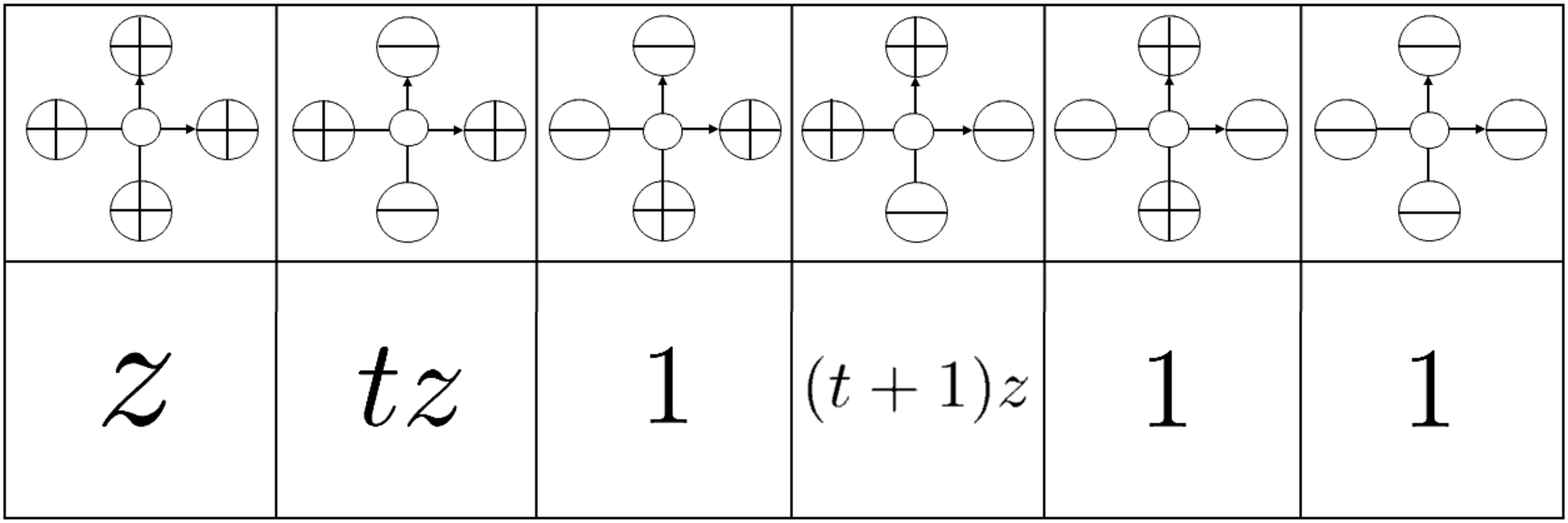}
\caption{The second $L$-operator \eqref{secondloperator}. 
We follow the pictorial description of \cite{Iv} for this second $L$-operator.}
\label{picturesecondloperator}
\end{figure}

\begin{figure}[ht]
\includegraphics[width=8cm]{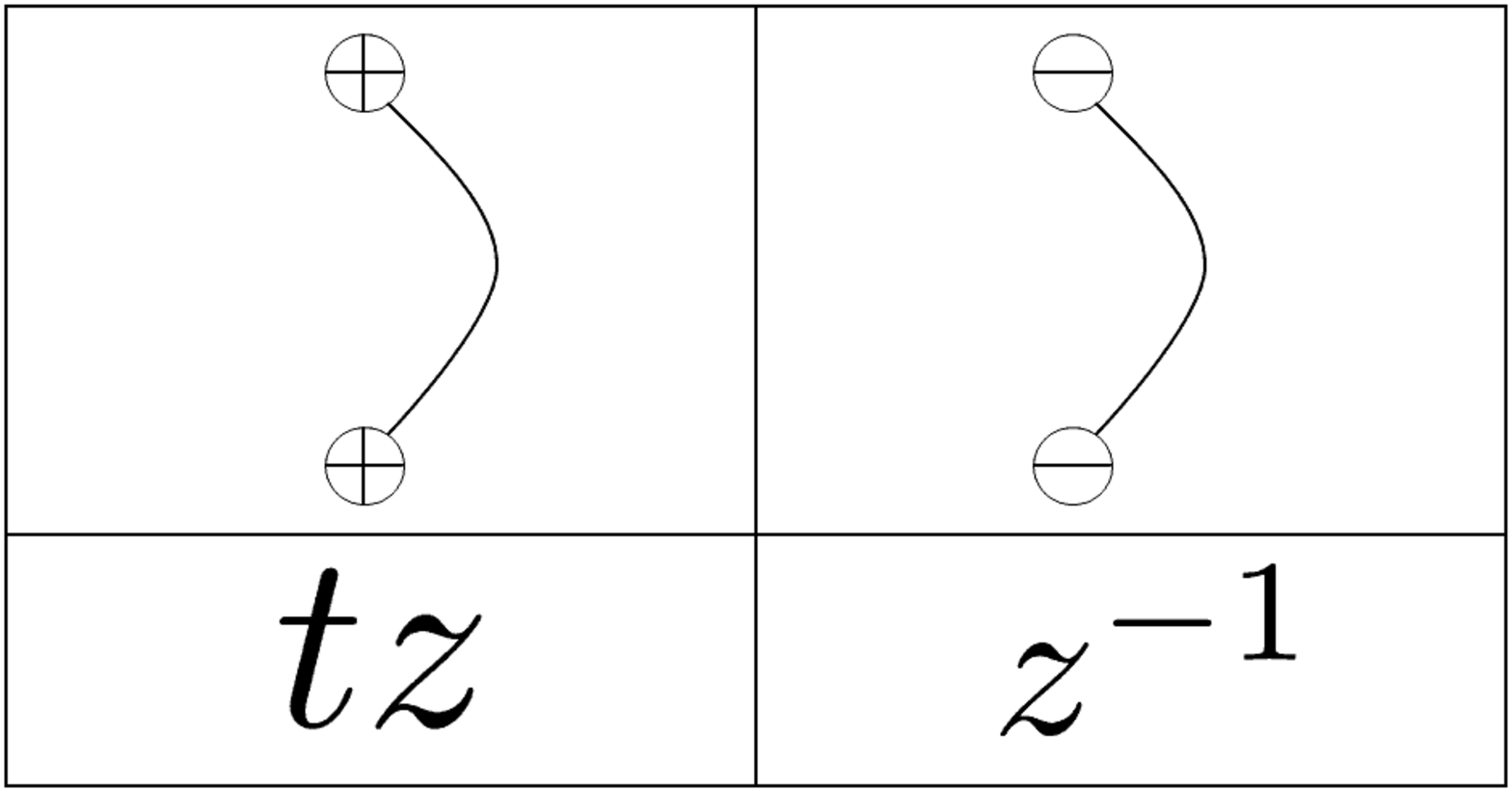}
\caption{The $K$-matrix \eqref{kmatrix}. 
We follow the pictorial description of \cite{Iv} for the $K$-matrix.}
\label{picturekmatrix}
\end{figure}

The $R$-matrix \eqref{rmatrix}
and $L$-operator \eqref{loperator} satsify the Yang-Baxter relation
\begin{align}
R_{ab}(z_1/z_2,t)L_{aj}(z_1,t)L_{bj}(z_2,t)
=L_{bj}(z_2,t)L_{aj}(z_1,t)R_{ab}(z_1/z_2,t), \label{RLL}
\end{align}
acting on $W_a \otimes W_b \otimes \mathcal{F}_j$.
We remark that this $RLL$ relation \eqref{RLL}
can be regarded as a special case of the
generalized Yang-Baxter relation for a more general
$R$-matrix \cite{Mu,DA,FCWZ}.
The $R$-matrix \eqref{rmatrix} and the $L$-operator
\eqref{loperator} in this section can be regarded as
different specializations of the general $R$-matrix
from this viewpoint.
One of the advantages of the point of view from the quantum group
is that one can systematically generalize the Felderhof model
to higher-dimensional representations \cite{DA}.

To realize the Tokuyama formula for Schur polynomials,
it was enough to use the $L$-operator
\eqref{loperator}.
To deal with symplectic Schur functions,
one needs more objects.
We introduce the second $L$-operator (see Figure
\ref{picturesecondloperator}
)
\begin{eqnarray}
\widetilde{L}_{aj}(z,t)=\left( 
\begin{array}{cccc}
z & 0 & 0 & 0 \\
0 & tz & 1 & 0 \\
0 & (t+1)z & 1 & 0 \\
0 & 0 & 0 & 1
\end{array}
\right), \label{secondloperator}
\end{eqnarray}
whose matrix elements are explicitly given by
\begin{align}
{}_a \langle 0| {}_j \langle 0 | \widetilde{L}_{a j}(z,t)
|0 \rangle_a | 0 \rangle_j&=z, \\
{}_a \langle 0| {}_j \langle 1 | \widetilde{L}_{a j}(z,t)
|0 \rangle_a | 1 \rangle_j&=tz, \\
{}_a \langle 0| {}_j \langle 1 | \widetilde{L}_{a j}(z,t)
|1 \rangle_a | 0 \rangle_j&=1, \\
{}_a \langle 1| {}_j \langle 0 | \widetilde{L}_{a j}(z,t)
|0 \rangle_a |1 \rangle_j&=(t+1)z, \\
{}_a \langle 1| {}_j \langle 0 | \widetilde{L}_{a j}(z,t)
|1 \rangle_a | 0 \rangle_j&=1, \\
{}_a \langle 1| {}_j \langle 1 | \widetilde{L}_{a j}(z,t)
|1 \rangle_a | 1 \rangle_j&=1.
\end{align}
We also introduce the $K$-matrix acting on the auxiliray space
$W_a$ (see Figure \ref{picturekmatrix})
\begin{eqnarray}
K_{a}(z,t)=\left( 
\begin{array}{cc}
tz & 0 \\
0 & z^{-1} \\
\end{array}
\right). \label{kmatrix}
\end{eqnarray}
The $K$-matrix is used when partition functions
of integrable lattice models under reflecting boundary conditions
are considered.
The matrix elements are explicitly given by
\begin{align}
{}_a \langle 0| K_{a}(z,t)
|0 \rangle_a&=tz, \\
{}_a \langle 1| K_{a}(z,t)
| 1 \rangle_a&=z^{-1}.
\end{align}

\begin{figure}[p]
\includegraphics[width=8cm]{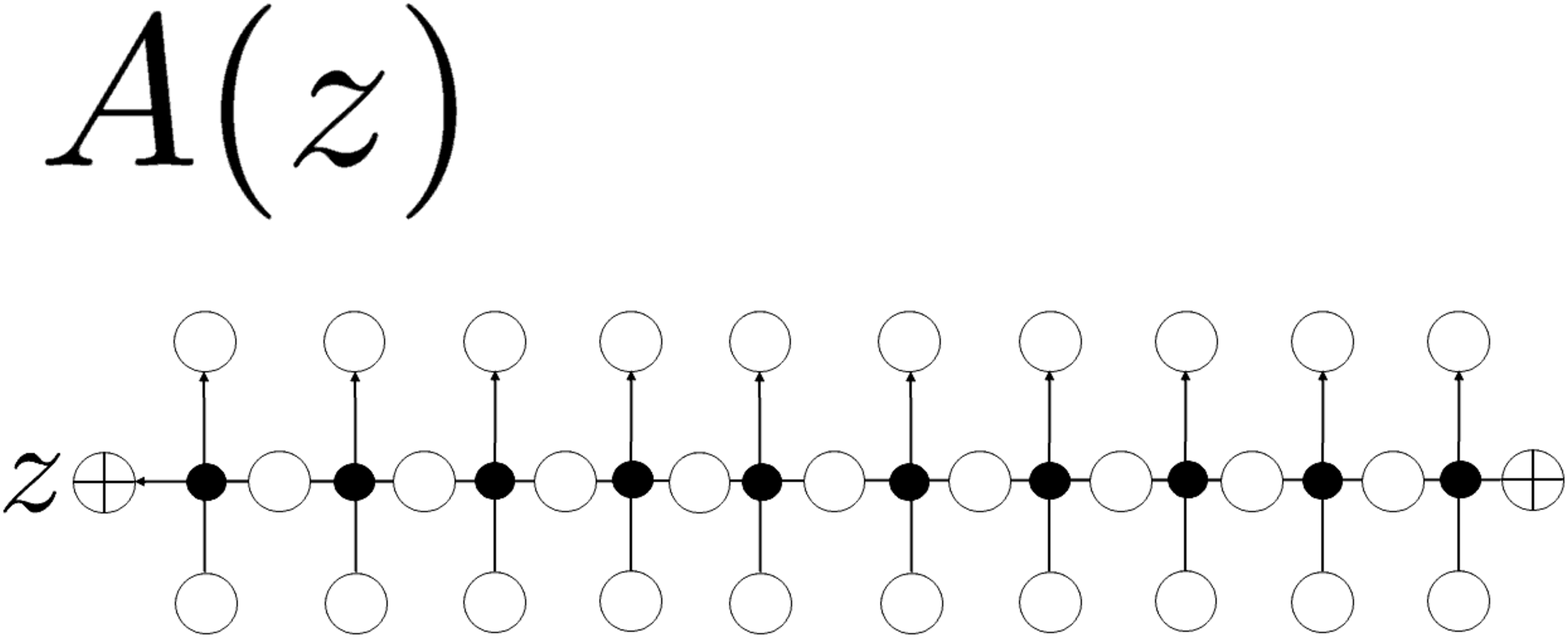}
\caption{The $A$-operator $A(z)$ 
\eqref{firstA}, which is a matrix element of
the monodromy matrix $T_a(z)$. The $A$-operator is
$2^M \times 2^M$ matrix-valued.
Both the leftmost and rightmost state on the horizontal line
(auxiliary space) are fixed as $\oplus$.}
\label{pictureAoperator}
\end{figure}
\begin{figure}[p]
\includegraphics[width=8cm]{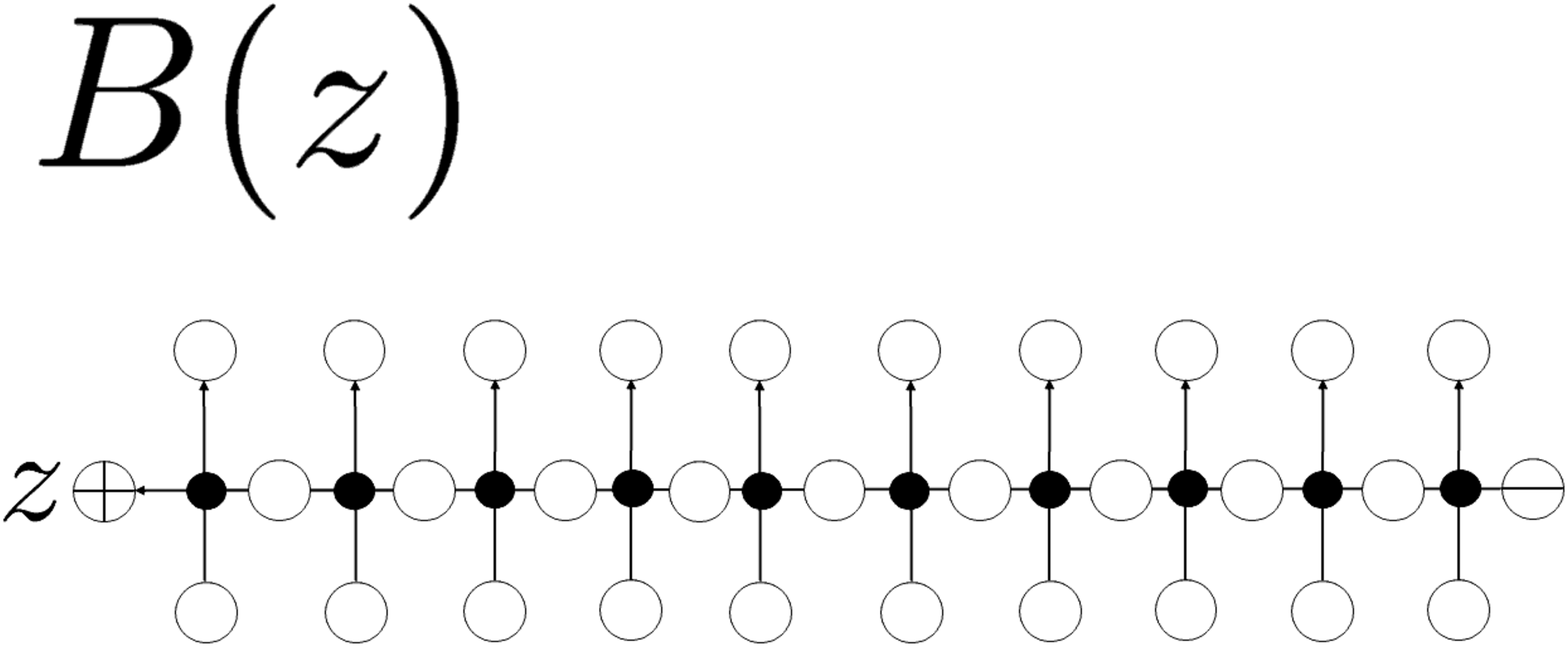}
\caption{The $B$-operator $B(z)$ \eqref{firstB}.
The leftmost state on the horizontal line
is fixed as $\oplus$,
whereas the rightmost state is fixed as $\ominus$.}
\label{pictureBoperator}
\end{figure}

\begin{figure}[p]
\includegraphics[width=8cm]{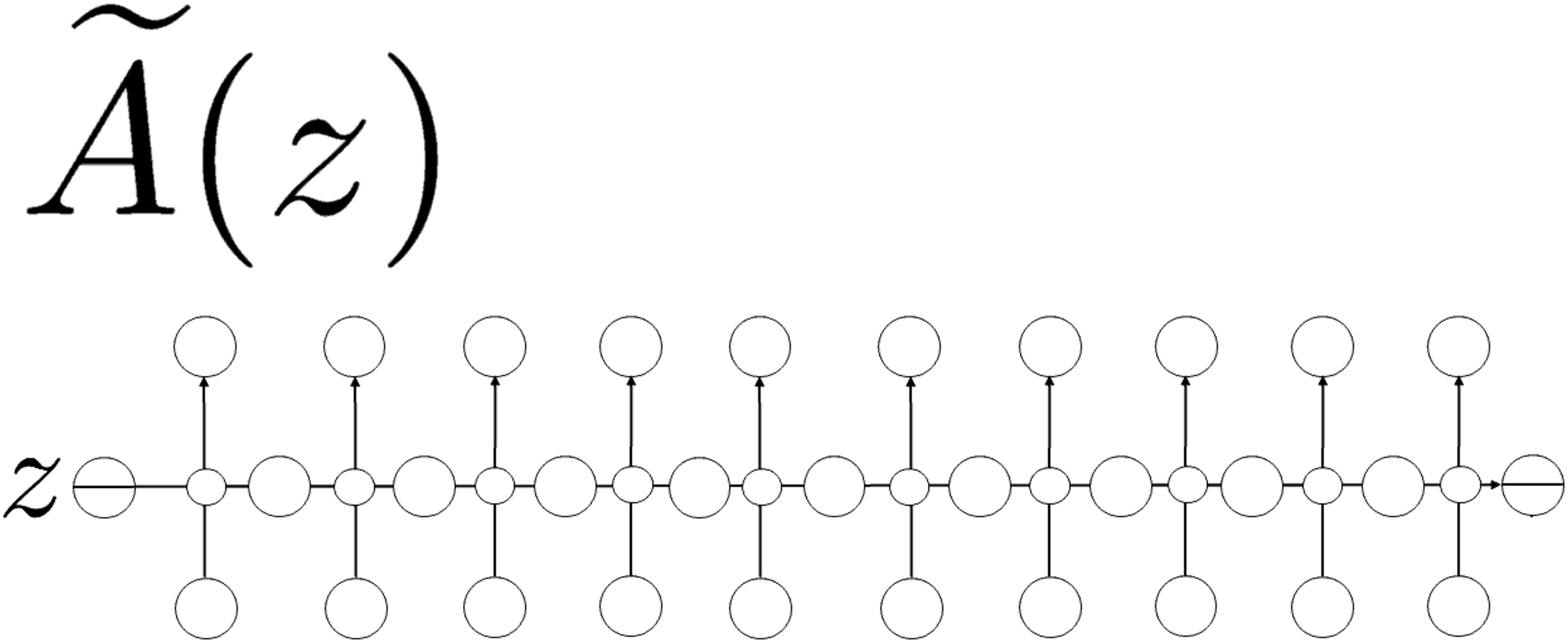}
\caption{The second $A$-operator $\widetilde{A}(z) \eqref{secondA}$.
Both the leftmost and rightmost state on the horizontal line
are fixed as $\ominus$.}
\label{picturesecondAoperator}
\end{figure}

\begin{figure}[ht]
\includegraphics[width=8cm]{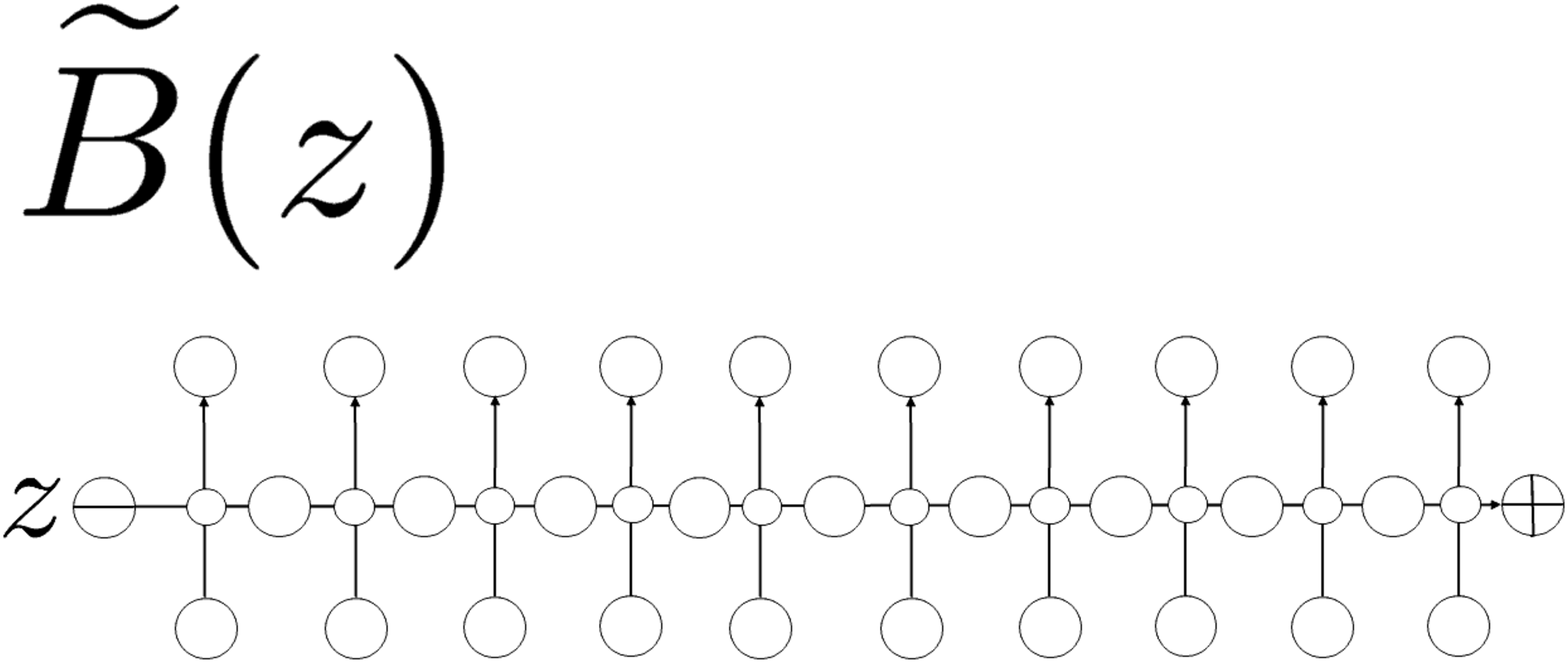}
\caption{The second $B$-operator $\widetilde{B}(z)$ \eqref{secondB}.
The leftmost state on the horizontal line
is fixed as $\ominus$,
whereas the rightmost state is fixed as $\oplus$.}
\label{picturesecondBoperator}
\end{figure}

From the $L$-operator, we construct two monodromy matrices
as products of $L$-operators
\begin{align}
T_{a}(z)=L_{a M}(z,t) \cdots L_{a 1}(z,t)
,
\label{monodromy1}
\end{align}
and
\begin{align}
\widetilde{T}_{a}(z)=\widetilde{L}_{a 1}(z,t) \cdots \widetilde{L}_{a M}(z,t)
,
\label{monodromy2}
\end{align}
which act on $W_a \otimes (\mathcal{F}_1\otimes\dots\otimes 
\mathcal{F}_{M})$.

In the language of quantum inverse scattering method,
the matrix elements of the monodromy matrices
with respect to the auxiliary space are called as $ABCD$ operators.
In this paper, we consider two types of
$A$- and $B$-operators
\begin{align}
A(z)={}_a \langle 0|T_{a}(z)|0 \rangle_{a}, \label{firstA} \\
B(z)={}_a \langle 0|T_{a}(z)|1 \rangle_{a}, \label{firstB}
\end{align}
and
\begin{align}
\widetilde{A}(z)={}_a \langle 1|\widetilde{T}_{a}(z)|1 \rangle_{a},
\label{secondA} \\
\widetilde{B}(z)={}_a \langle 0|\widetilde{T}_{a}(z)|1 \rangle_{a}.
\label{secondB}
\end{align}

See Figures \ref{pictureAoperator}, \ref{pictureBoperator},
\ref{picturesecondAoperator} and \ref{picturesecondBoperator}
for a graphical description of these two types of $A$- and $B$-operators.

\begin{figure}[ht]
\includegraphics[width=12cm]{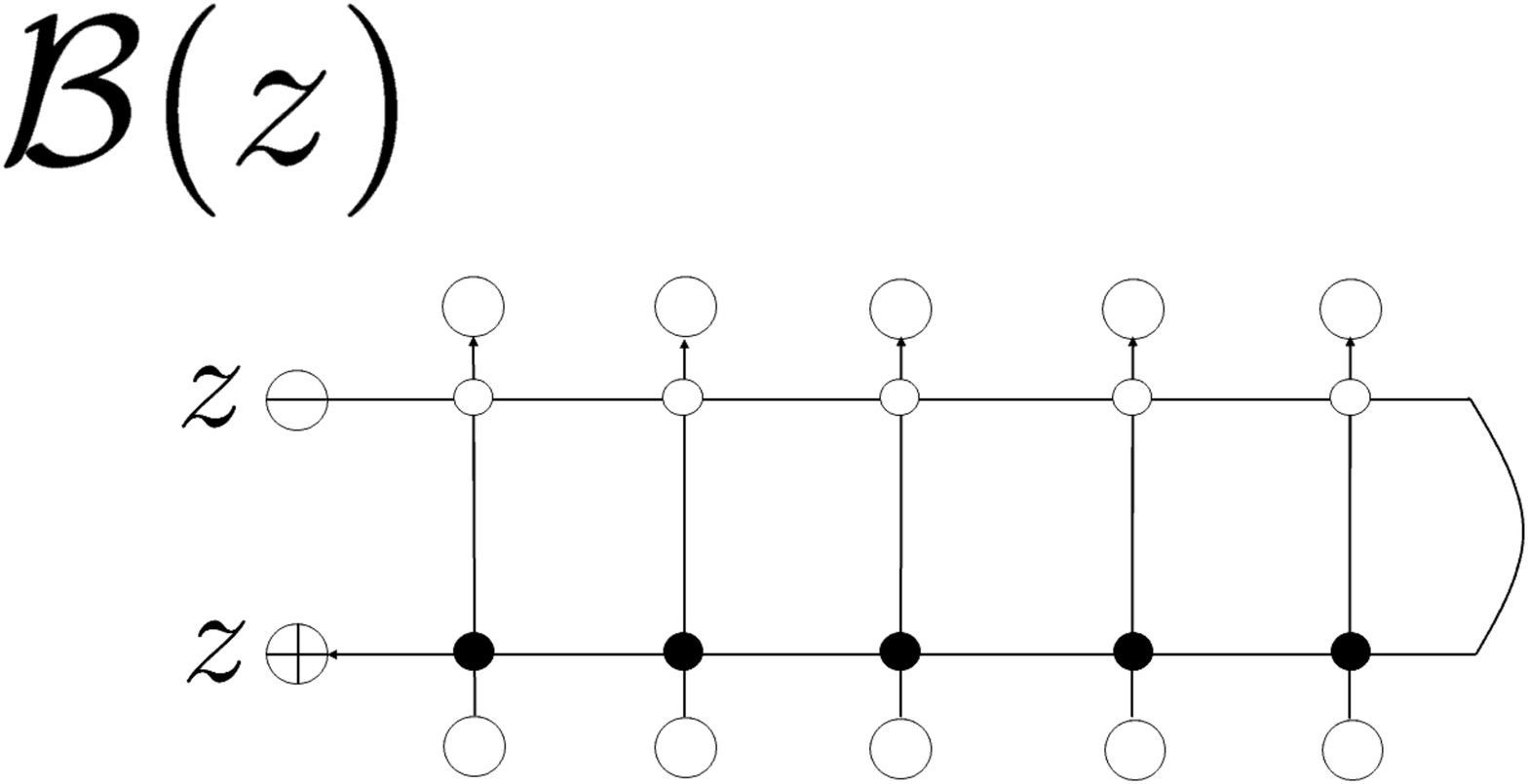}
\caption{The double row $B$-operator $\mathcal{B}(z)$ (the left hand side of
\eqref{doublerow}).}
\label{picturedoublerowBoperator}
\end{figure}

\begin{figure}[ht]
\includegraphics[width=12cm]{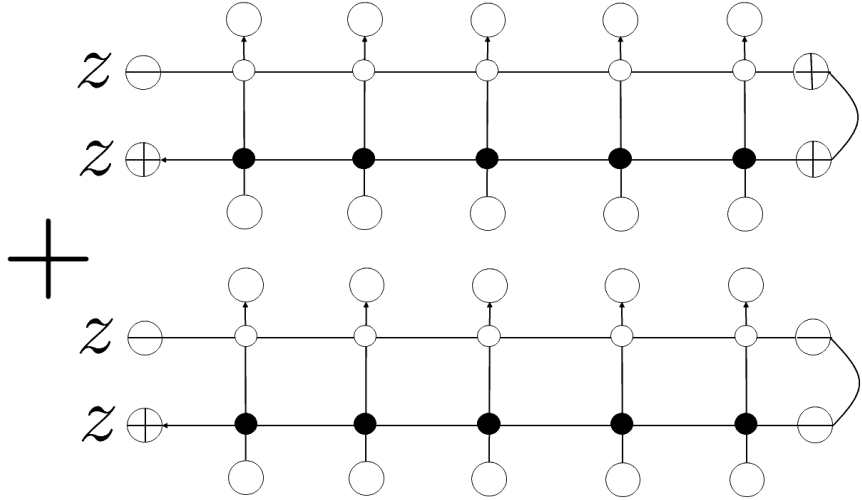}
\caption{The decomposition of the double row $B$-operator $\mathcal{B}(z)$ 
(the right hand side of \eqref{doublerow}).
The top figure where the rightmost states are fixed as $\oplus$
represents the term $\widetilde{B}(z) {}_a \langle 0| K(z,t)|0 \rangle_a A(z)
=tz \widetilde{B}(z)A(z)$.
The bottom figure where the rightmost states are fixed as $\oplus$
represents the term $\widetilde{A}(z) {}_a \langle 1| K(z,t)|1 \rangle_a B(z)
=z^{-1} \widetilde{A}(z) B(z)$.
}
\label{picturedoublerowBoperatordecomposition}
\end{figure}

By using these four monodromy operators and $K$-matrix,
we introduce the following double row $B$-operator \cite{Sklyanin}
which we use to construct the wavefunctions
under reflecting boundary.

\begin{align}
\mathcal{B}(z)&=\widetilde{B}(z) {}_a \langle 0| K(z,t)|0 \rangle_a A(z)
+\widetilde{A}(z) {}_a \langle 1| K(z,t)|1 \rangle_a B(z) \\
&=tz \widetilde{B}(z)A(z)
+z^{-1} \widetilde{A}(z) B(z). \label{doublerow}
\end{align}

See Figures \ref{picturedoublerowBoperator} and
\ref{picturedoublerowBoperatordecomposition}
for pictorial descriptions of \eqref{doublerow}.

\section{Wavefunction and symplectic Schur functions}
We introduce the wavefunction which is a special class of partition functions,
and review the relation with the symplectic Schur functions defined below.
\begin{definition}
The symplectic Schur functions is defined to be the following determinant:
\begin{align}
sp_\lambda(\{ z \}_N)=
   \frac{\mathrm{det}_N(z_j^{\lambda_k+N-k+1}-z_j^{-\lambda_k-N+k-1})}
        {\mathrm{det}_N(z_j^{N-k+1}-z_j^{-N+k-1})},
 \label{symplecticSchur}
\end{align}
where $\{ z \}_N=\{z_1,\dots,z_N \}$ is a set of variables
and $\lambda$ denotes a Young diagram
$\lambda=(\lambda_1,\lambda_2,\dots,\lambda_N)$
with weakly decreasing non-negative integers
$\lambda_1 \ge \lambda_2 \ge \cdots \ge \lambda_N \ge 0$.
\end{definition}
Before introducing the wavefunction, we first define
the arbitrary $N$-particle state \\
$|\Psi(z_1,\dots,z_N) \ket$
 with $N$ spectral parameters
$\{ z \}_N=\{ z_1,\dots,z_N \}$
as a state constructed by a multiple action
of the double row $B$-operator on the vacuum state 
$|\Omega \ket:= | 0^{M} \ket:=|0\ket_1\
\otimes \dots \otimes |0\ket_{M}$
\begin{align}
|\Psi(z_1,\dots,z_N) \ket=\mathcal{B}(z_1)
\cdots \mathcal{B}(z_N)| \Omega \ket.
\label{statevector}
\end{align}

\begin{figure}[ht]
\includegraphics[width=15cm]{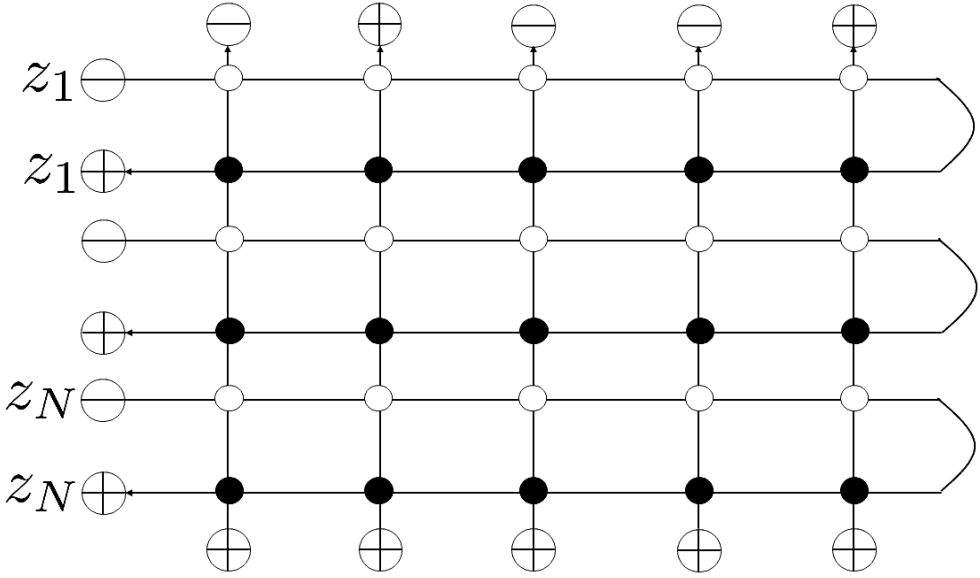}
\caption{The wavefunction
$\langle x_1 \cdots x_N | \Psi(z_1,\dots,z_N) \ket$
for the case $M=5$, $N=3$, $(x_1,x_2,x_3)
=(2,3,5)$.}
\label{picturewavefunction}
\end{figure}

Next, we introduce the wavefunction
$\bra x_1 \cdots x_N | \Psi(z_1,\dots,z_N) \ket$
as the overlap between an arbitrary off-shell 
$N$-particle state $|\Psi(z_1,\dots,z_N)\ket$ and 
the (normalized) state with an arbitrary particle configuration 
$|x_1 \cdots x_N\ket$ $(1 \le x_1<\dots<x_N \le M$), 
where $x_j$ denotes the positions of the particles. 
The particle configurations are explicitly defined as
\begin{align}
   \bra x_1 \cdots x_N|=\bra \Omega|\prod_{j=1}^N \sigma^+_{x_j},
\end{align}
where $\langle \Omega|:=\langle 0^{M}|:=
{}_1\bra 0|\otimes\dots \otimes{}_M\bra 0|$.
Here, we define $\sigma^+$ and $\sigma^-$ as
operators acting on the basis elements as
\begin{align}
&\sigma^+|1 \rangle=|0 \rangle, \ 
\sigma^+|0 \rangle=0, \ 
\langle 0|\sigma^+=\langle 1|, \
\langle 1|\sigma^+=0, 
\\
&\sigma^-|0 \rangle=|1 \rangle, \
\sigma^-|1 \rangle=0, \
\langle 1|\sigma^-=\langle 0|, \
\langle 0|\sigma^-=0.
\end{align}
The subscript $j$ of $\sigma_j^+$ or $\sigma_j^-$
indicates  that the operator acts on the space $\mathcal{F}_j$ as $\sigma^+$
or $\sigma^-$,
and as an idenitity on the other spaces.

The following correspondence
between the wavefunction of the Felderhof model
under reflecting boundary
and the symplectic Schur functions was proved by Ivanov.
\begin{theorem} \label{TheoremIvanov} \cite{Iv}
The wavefunction $\bra x_1 \cdots x_N | \Psi(z_1,\dots,z_N) \ket$
is expressed by the symplectic Schur functions as
\begin{align}
\bra x_1 \cdots x_N | \Psi(z_1,\dots,z_N) \ket
=\prod_{j=1}^N z_j^{j-1-N}(1+tz_j^2)
\prod_{1 \le j<k \le N}(1+tz_j z_k)(1+tz_j z_k^{-1})
sp_\lambda(\{ z \}_N).
\end{align}
Here the Young diagram for the symplectic Schur functions correspond to
the particle configuration under the relation
$\lambda_j=x_{N-j+1}-N+j-1$, $j=1,\dots,N$.
\end{theorem}
The above theorem means that the product of a deformation of Weyl's denominator
and the symplectic Schur functions, which is
an irreducible character of the symplectic group $Sp(2n,\mathbf{C})$,
can be expressed as a wavefunction of the free-fermion six-vertex model
under the reflecting boundary condition.
The wavefunction offers an explicit description
in terms the Proctor pattern, and hence
this result can be regarded as a type $C$ version of the Tokuyama formula.

\section{Dual wavefunction}
We now introduce the dual wavefunction, and
study the exact relation between it
and the symplectic Schur functions.
The strategy of the proof of the correspondence given here
can be regarded as the symplectic version of our proof
of the correspondence between the dual wavefunction
without reflecting boundary
and the Schur polynomials \cite{LMP}.
We transform the statement of the theorem into another equivalent one
which enables us to use the arguments of Bump-Brubaker-Friedberg \cite{BBF}
and Ivanov \cite{Iv}.

Before defining the dual wavefunction,
we introduce another type of arbitrary dual $N$-hole state
$\langle \Phi(z_1,\dots,z_N)|$ by a multiple action of
the double row $B$-operator on the dual particle occupied state
$\langle 1 \cdots M|:=\langle 1^{M}|:=
{}_1\bra 1|\otimes\dots \otimes{}_M\bra 1|$
\begin{align}
\langle \Phi(z_1,\dots,z_N)|
=\langle 1 \cdots M|\mathcal{B}(z_1) \cdots \mathcal{B}(z_N).
\end{align}
It is convenient to introduce a notation for
the state with an arbitrary hole configuration 
$|\overline{x_1} \cdots \overline{x_N} \ket$
$(1 \le \overline{x_1}<\dots< \overline{x_N} \le M$), 
where $\overline{x_j}$ denotes the positions of holes. 
Explicitly,
\begin{align}
|\overline{x_1} \cdots \overline{x_N} \ket
=\prod_{j=1}^N \sigma^+_{x_j}
(|1 \rangle_1 \otimes \cdots \otimes |1 \rangle_M).
\end{align}

\begin{figure}[ht]
\includegraphics[width=15cm]{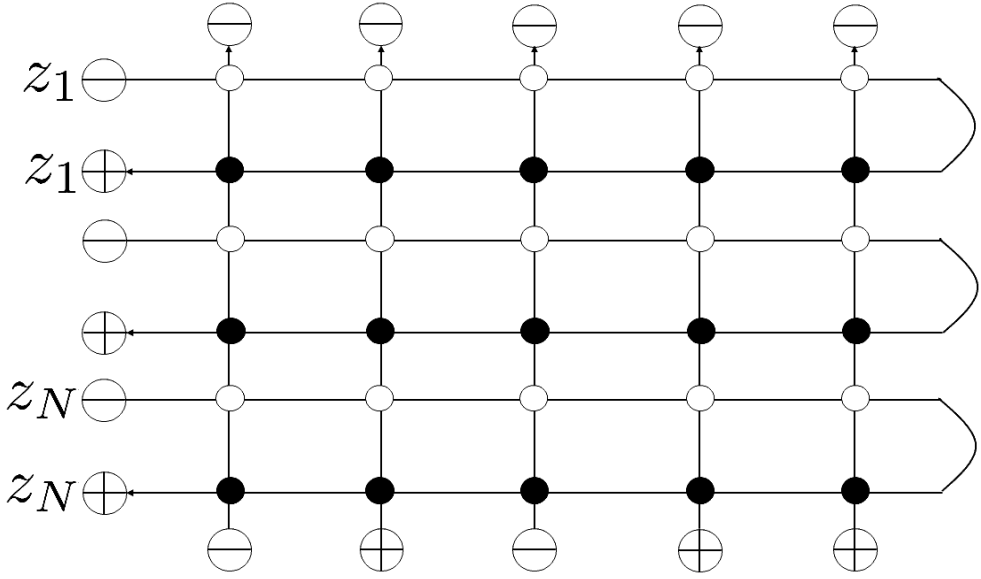}
\caption{The dual wavefunction
$\langle \Phi(z_1,\dots,z_N)|\overline{x_1} \cdots \overline{x_N} \ket$
for the case $M=5$, $N=3$, $(\overline{x_1},\overline{x_2},\overline{x_3})
=(1,2,4)$.}
\label{picturedualwavefunction}
\end{figure}

The dual wavefunction
$\langle \Phi(z_1,\dots,z_N)|\overline{x_1} \cdots \overline{x_N} \ket$
is defined as the overlap between the
arbitrary dual $N$-hole state
$\langle \Phi(z_1,\dots,z_N)|$
and hole configurations
$|\overline{x_1} \cdots \overline{x_N} \ket$
(see Figure \ref{picturedualwavefunction} for an example
of a graphical description of the dual wavefunction).

We show the following relation between the dual wavefunction
and the symplectic Schur functions.
\begin{theorem} \label{theoremdualandsymplecticschur}
The dual wavefunction
$\langle \Phi(z_1,\dots,z_N)| \overline{x_1} \cdots \overline{x_N}
\rangle$ can be expressed by the symplectic Schur functions as
\begin{align}
&\langle \Phi(z_1,\dots,z_N)| \overline{x_1} \cdots \overline{x_N}
\rangle \nonumber \\
=&t^{N(M-N)}
\prod_{j=1}^N z_j^{j-1-N}(1+tz_j^2)
\prod_{1 \le j<k \le N}(1+tz_j z_k)(1+tz_j z_k^{-1})
sp_{\overline{\lambda}} ( \{ tz \}_N ).
\label{dualwavefunctionsandsymplecticschur}
\end{align}
Here the Young diagram for the symplectic Schur functions correspond to
the hole configuration under the relation
$\overline{\lambda_j}=\overline{x_{N-j+1}}-N+j-1$, $j=1,\dots,N$,
and the symmetric variables are
$\displaystyle \{ tz \}_N=
\{ tz_1,\dots,tz_N \}$.
\end{theorem}
The result of Theorem \ref{theoremdualandsymplecticschur} resembles
the case for the dual wavefunction without reflecting boundary
which gives the Schur polynomials.
There is again a factor $t^{N(M-N)}$ which depends on the number
of sites $M$ and the number of particles $N$ in the right hand side
of \eqref{dualwavefunctionsandsymplecticschur}.
Also, the symmetric variables of the symplectic Schur functions
are $\displaystyle \{ tz \}_N$, not simply $\displaystyle \{ z \}_N$.

It seems difficult to directly prove
\eqref{dualwavefunctionsandsymplecticschur} itself by using the argument by
Ivanov \cite{Iv}, which extends the one for the case without reflecting boundary
by Bump-Brubaker-Friedberg \cite{BBF} to the symplectic ice.
As in the case of the paper which we proved the correspondence
between the dual wavefunction without reflecting boundary
and the Schur polynomials,
the key is to transform the statement 
\eqref{dualwavefunctionsandsymplecticschur} to another equivalent one,
which enables us to give a proof by using the argument by Ivanov.
\begin{proof}
We transform \eqref{dualwavefunctionsandsymplecticschur}
as follows. First,
by rescaling each $z_j$ to $t^{-1} z_j$,
we have
\begin{align}
&\langle \Phi(t^{-1}z_1,\dots,t^{-1}z_N)
| \overline{x_1} \cdots \overline{x_N}
\rangle \nonumber \\
=&t^{MN}
\prod_{j=1}^N z_j^{j-1-N}(1+t^{-1}z_j^2)
\prod_{1 \le j<k \le N}(1+t^{-1}z_j z_k)(t^{-1}+z_j z_k^{-1})
sp_{\overline{\lambda}} ( \{ z \}_N ).
\label{dualwavefunctionsandsymplecticschurtwo}
\end{align}
We further rewrite
\eqref{dualwavefunctionsandsymplecticschurtwo}
in the following form.
\begin{align}
&t^N
\langle 1 \cdots M|\Bigg(\frac{\mathcal{B}(t^{-1}z_1)}{t^{M+1}}\Bigg) \cdots
\Bigg(\frac{\mathcal{B}(t^{-1}z_N)}{t^{M+1}} \Bigg)| \overline{x_1} \cdots \overline{x_N}
\rangle \nonumber \\
=&
\prod_{j=1}^N z_j^{j-1-N}(1+t^{-1}z_j^2)
\prod_{1 \le j<k \le N}(1+t^{-1}z_j z_k)(t^{-1}+z_j z_k^{-1})
sp_{\overline{\lambda}} ( \{ z \}_N ).
\label{equivalentone}
\end{align}
For giving a proof, it is convenient to
introduce the following rescaled $L$-operators and $K$-matrix
\begin{eqnarray}
L^\prime(z,t)=
\frac{1}{t}L(t^{-1}z,t)=
\left( 
\begin{array}{cccc}
t^{-1} & 0 & 0 & 0 \\
0 & 1 & t^{-1} & 0 \\
0 & (t+1)z^{-1} & z^{-1} & 0 \\
0 & 0 & 0 & z^{-1}
\end{array}
\right), \label{loperator2}
\end{eqnarray}

\begin{eqnarray}
\widetilde{L}^\prime(z,t)=
\widetilde{L}(t^{-1}z,t)=
\left( 
\begin{array}{cccc}
t^{-1}z & 0 & 0 & 0 \\
0 & z & 1 & 0 \\
0 & (1+t^{-1})z & 1 & 0 \\
0 & 0 & 0 & 1
\end{array}
\right), \label{loperator3}
\end{eqnarray}

\begin{eqnarray}
K^\prime(z,t)
=\frac{1}{t}K(t^{-1}z)=\left( 
\begin{array}{cc}
t^{-1}z & 0 \\
0 & z^{-1} \\
\end{array}
\right). \label{kmatrix2}
\end{eqnarray}

We denote the two types of the $A$- and $B$-operators,
the double row $B$-operator
constructed from these rescaled $L$-operators and $K$-matrix
$L^\prime(z,t)$, $\widetilde{L}^\prime(z,t)$, $K^\prime(z,t)$
as $A^\prime(z)$, $B^\prime(z)$, $\widetilde{A}^\prime(z)$,
$\widetilde{B}^\prime(z)$ and $\mathcal{B}^\prime(z)$ respectively.
The double row $B$ operator is now constructed from the
two types of the $A$- and $B$-operators as
\begin{align}
\mathcal{B}^\prime(z)&=\widetilde{B}^\prime(z) {}_a \langle 0| K^\prime(z,t)|0 \rangle_a A^\prime(z)
+\widetilde{A}^\prime(z) {}_a \langle 1| K^\prime(z,t)|1 \rangle_a B^\prime(z) \\
&=t^{-1}z \widetilde{B}^\prime(z)A^\prime(z)
+z^{-1} \widetilde{A}^\prime(z) B^\prime(z). \label{doublerow2}
\end{align}

Using these rescaled objects,
\eqref{equivalentone} can be expressed as
\begin{align}
&t^N
\langle 1 \cdots M|\mathcal{B}^\prime(z_1) \cdots
\mathcal{B}^\prime(z_N)| \overline{x_1} \cdots \overline{x_N}
\rangle \nonumber \\
=&
\prod_{j=1}^N z_j^{j-1-N}(1+t^{-1}z_j^2)
\prod_{1 \le j<k \le N}(1+t^{-1}z_j z_k)(t^{-1}+z_j z_k^{-1})
sp_{\overline{\lambda}} ( \{ z \}_N )
. \label{equivalenttwo}
\end{align}
Instead of proving \eqref{dualwavefunctionsandsymplecticschur},
we show \eqref{equivalenttwo} since this is equivalent to
\eqref{dualwavefunctionsandsymplecticschur} and is the expression
which one can use the argument given in Ivanov \cite{Iv}.

We first show the following lemma.

\begin{lemma} \label{lemmaone}
\begin{align}
&\prod_{j=1}^N z_j^{N+1-j} (1+t^{-1}z_j^2)^{-1}
\prod_{1 \le j<k \le N}(1+t^{-1}z_j z_k)^{-1}
(t^{-1}+z_j z_k^{-1})^{-1} \nonumber \\
\times&t^N
\langle 1 \cdots M|\mathcal{B}^\prime(z_1) \cdots
\mathcal{B}^\prime(z_N)| \overline{x_1} \cdots \overline{x_N}
\rangle
\end{align}
does not depend on $t$.
\end{lemma}
\begin{proof}
We follow along the lines of the proof by Ivanov \cite{Iv}.
It is enough to prove the following properties
for \\
$
t^N
\langle 1 \cdots M|\mathcal{B}^\prime(z_1) \cdots
\mathcal{B}^\prime(z_N)| \overline{x_1} \cdots \overline{x_N}
\rangle
$: 
\begin{enumerate}
\item 
$
t^N
\langle 1 \cdots M|\mathcal{B}^\prime(z_1) \cdots
\mathcal{B}^\prime(z_N)| \overline{x_1} \cdots \overline{x_N}
\rangle
$
is a polynomial of $t^\prime:=t^{-1}$ with its highest degree at most $N^2$. 
\item
$
t^N
\langle 1 \cdots M|\mathcal{B}^\prime(z_1) \cdots
\mathcal{B}^\prime(z_N)| \overline{x_1} \cdots \overline{x_N}
\rangle/D^\prime
$, \\
where $D^\prime:=\prod_{j=1}^N z_j^{j-1-N} (1+t^\prime z_j^2)
\prod_{1 \le j < k \le N}(1+t^\prime z_j z_k)
(t^{\prime}+z_j z_k^{-1})$, is invariant under
any permutation of $z_1,\dots,z_N$.
\item
$
t^N
\langle 1 \cdots M|\mathcal{B}^\prime(z_1) \cdots
\mathcal{B}^\prime(z_N)| \overline{x_1} \cdots \overline{x_N}
\rangle/D^\prime
$ is invariant under $z_i \longleftrightarrow z_i^{-1}$.
\end{enumerate}
We first show
$
\mathrm{deg}_{t^\prime}
(t^N
\langle 1 \cdots M|\mathcal{B}^\prime(z_1) \cdots
\mathcal{B}^\prime(z_N)| \overline{x_1} \cdots \overline{x_N}
\rangle)
\le N^2
$ by induction on $N$.
We use the following properties for the matrix elements
of the $A$- and $B$-operators,
which can easily be seen from the matrix elements 
of the $L$-operators.
\begin{align}
&0 \le \mathrm{deg}_{t^\prime}(t
\langle \overline{x_1} \cdots \overline{x_N} |A^\prime(z)|
\overline{y_1}
\cdots \overline{y_{N}} \rangle
) \le N-1, \ \ \ N \ge 1, \label{degpropone} \\
&0 \le \mathrm{deg}_{t^\prime}(t
\langle \overline{x_1} \cdots \overline{x_N}
|B^\prime(z)|
\overline{y_1}
\cdots \overline{y_{N+1}} \rangle
) \le N, \ \ \ N \ge 0, \label{degproptwo} \\
&0 \le \mathrm{deg}_{t^\prime}(
\langle \overline{x_1} \cdots \overline{x_N} |\widetilde{A}^\prime(z)|
\overline{y_1}
\cdots \overline{y_{N}} \rangle
) \le N, \ \ \ N \ge 0, \label{degpropthree} \\
&0 \le \mathrm{deg}_{t^\prime}(
\langle \overline{x_1} \cdots \overline{x_N}
|\widetilde{B}^\prime(z)|
\overline{y_1}
\cdots \overline{y_{N+1}} \rangle
) \le N, \ \ \ N \ge 0. \label{degpropfour}
\end{align}
We first remark that it is enough to consider the above matrix elements,
since all of the non-zero matrix elements are included in the above cases.
For example, one does not need to consider
$\langle \overline{x_1} \cdots \overline{x_N} |A^\prime(z)|
\overline{y_1}
\cdots \overline{y_{N+1}} \rangle$.
This is 
because due to the so-called ice rule 
${}_a \langle \delta| {}_j \langle \gamma|L(z)|\alpha \rangle_a |\beta \rangle_j=0$
unless $\alpha+\beta=\gamma+\delta$
for the $L$-operator of the six-vertex model,
the $A$-operators preserve the total number of holes in the quantum space,
hence \\
$\langle \overline{x_1} \cdots \overline{x_N} |A^\prime(z)|
\overline{y_1} \cdots \overline{y_{N+1}} \rangle=0$ follows immediately.
Similary, the $B$-operators increase the total number of holes
in the quantum space by one, one can see
$\langle \overline{x_1} \cdots \overline{x_N} |B^\prime(z)|
\overline{y_1} \cdots \overline{y_{N}} \rangle=0$ immediately.
It is not difficult to find the above properties of the degree
with respect to $t^\prime$ by taking the ice-rule into account.

Let us show Property 1 for the case $N=1$.
We use the decomposition of the double row $B$-operator \eqref{doublerow2}
and insert the completeness relation
$\displaystyle \sum_{\overline{y}}|\overline{y} \rangle \langle \overline{y}|=
\mathrm{Id}$ between the $A$- and $B$-operators
to deform $\mathrm{det}_{t^\prime}(t
\langle 1 \cdots M|\mathcal{B}^\prime(z)|
 \overline{x}
\rangle)$ as
\begin{align}
&\mathrm{det}_{t^\prime}(t
\langle 1 \cdots M|\mathcal{B}^\prime(z)|
 \overline{x}
\rangle)
=\mathrm{det}_{t^\prime}(t
\langle 1 \cdots M|
(t^{-1}z \widetilde{B}^\prime(z)A^\prime(z)
+z^{-1} \widetilde{A}^\prime(z) B^\prime(z))
|
\overline{x}
\rangle) \nonumber \\
=&\mathrm{deg}_{t^\prime}
(t^{-1}z \sum_{\overline{y}}
\langle 1 \cdots M| \widetilde{B}^\prime(z) | \overline{y} \rangle
\langle \overline{y} | tA^\prime(z)|\overline{x} \rangle
+z^{-1}
\langle 1 \cdots M| \widetilde{A}^\prime(z) | 1 \cdots M \rangle
\langle 1 \cdots M| tB^\prime(z)|\overline{x} \rangle
).
\end{align}
Using the generic property of the degree
$\mathrm{deg}(P+Q)=\mathrm{Max}(\mathrm{deg}P,\mathrm{deg}Q)$
and the properties
\eqref{degpropone}, \eqref{degproptwo}, \eqref{degpropthree}
and \eqref{degpropfour},
it follows that $0 \le \mathrm{deg}_{t^\prime}(t
\langle 1 \cdots M|\mathcal{B}^\prime(z)|
 \overline{x}
\rangle) \le 1$.

Next, let us assume
$
\mathrm{deg}_{t^\prime}
(t^N
\langle 1 \cdots M|\mathcal{B}^\prime(z_1) \cdots
\mathcal{B}^\prime(z_N)| \overline{x_1} \cdots \overline{x_N}
\rangle)
\le N^2
$.
Let us show
$
0 \le \mathrm{deg}_{t^\prime}
(t \langle \overline{x_1} \cdots \overline{x_N} |\mathcal{B}^\prime(z)| \overline{y_1}
\cdots \overline{y_{N+1}} \rangle) \le 2N+1
$.
This can be seen as above by using the 
decomposition of the double row $B$-operator \eqref{doublerow2}
and inserting the completeness relation
to deform $\mathrm{deg}_{t^\prime}
(t \langle \overline{x_1} \cdots \overline{x_N} |\mathcal{B}^\prime(z)| \overline{y_1}
\cdots \overline{y_{N+1}} \rangle)$
as
\begin{align}
&\mathrm{deg}_{t^\prime}
(t \langle \overline{x_1} \cdots \overline{x_N} |\mathcal{B}^\prime(z)| \overline{y_1}
\cdots \overline{y_{N+1}} \rangle) \nonumber \\
=&\mathrm{deg}_{t^\prime}
(t \langle \overline{x_1} \cdots \overline{x_N} |
(t^{-1}z \widetilde{B}^\prime(z)A^\prime(z)
+z^{-1} \widetilde{A}^\prime(z) B^\prime(z))
| \overline{y_1}
\cdots \overline{y_{N+1}} \rangle) \nonumber \\
=&\mathrm{deg}_{t^\prime}\Bigg(
t^{-1}z \sum_{\{ \overline{u} \}}
\langle \overline{x_1} \cdots \overline{x_N}|
\widetilde{B}^\prime(z)| \overline{u_1} \cdots \overline{u_{N+1}}
\rangle
\langle \overline{u_1} \cdots \overline{u_{N+1}}|t A^\prime(z)|
\overline{y_1}
\cdots \overline{y_{N+1}} \rangle \nonumber \\
+&z^{-1}
\sum_{\{ \overline{v} \}}
\langle \overline{x_1} \cdots \overline{x_N}|
\widetilde{A}^\prime(z)| \overline{v_1} \cdots \overline{v_{N}}
\rangle
\langle \overline{v_1} \cdots \overline{v_{N}}|t B^\prime(z)|
\overline{y_1}
\cdots \overline{y_{N+1}} \rangle
\Bigg),
\end{align}
and using $\mathrm{deg}(P+Q)=\mathrm{Max}(\mathrm{deg}P,\mathrm{deg}Q)$,
the properties \eqref{degpropone},
\eqref{degproptwo}, \eqref{degpropthree} and \eqref{degpropfour}.

One can finally see
$
\mathrm{deg}_{t^\prime}
(t^{N+1}
\langle 1 \cdots M|\mathcal{B}^\prime(z_1) \cdots
\mathcal{B}^\prime(z_{N+1})| \overline{y_1} \cdots \overline{y_{N+1}}
\rangle)
\le (N+1)^2
$
by using the property
$
0 \le \mathrm{deg}_{t^\prime}
(t \langle \overline{x_1} \cdots \overline{x_N} |\mathcal{B}^\prime(z)| \overline{y_1}
\cdots \overline{y_{N+1}} \rangle) \le 2N+1
$
together with the assumption
$
\mathrm{deg}_{t^\prime}
(t^N
\langle 1 \cdots M|\mathcal{B}^\prime(z_1) \cdots
\mathcal{B}^\prime(z_N)| \overline{x_1} \cdots \overline{x_N}
\rangle)
\le N^2
$
and the decomposition
\begin{align}
&t^{N+1}
\langle 1 \cdots M|\mathcal{B}^\prime(z_1) \cdots
\mathcal{B}^\prime(z_{N+1})| \overline{y_1} \cdots \overline{y_{N+1}}
\rangle \nonumber \\
=&\sum_{\{ \overline{x} \}}
(t^N
\langle 1 \cdots M|\mathcal{B}^\prime(z_1) \cdots
\mathcal{B}^\prime(z_{N})| \overline{x_1} \cdots \overline{x_{N}} \rangle)
(t \langle \overline{x_1} \cdots \overline{x_N} |\mathcal{B}^\prime(z_{N+1})| \overline{y_1}
\cdots \overline{y_{N+1}} \rangle).
\end{align}

Property 2 can be proved exactly
in the same way as in Lemma 2 in Ivanov \cite{Iv} which is a lengthy one,
since one needs many local relations and arguments to prove the lemma.
Instead, given below is a shortcut of the argument by
transforming his result to Property 2.
First, one notes that Ivanov's proof of Lemma 2 shows that
not only $\langle x_1 \cdots x_N|\mathcal{B}(z_1) \cdots
\mathcal{B}(z_N)|\Omega \rangle$
but also every matrix element
$ \langle y_1 \cdots y_{n+N} |\mathcal{B}(z_1) \cdots \mathcal{B}(z_N) 
|x_1 \cdots x_n \rangle$ has the property that
$
\langle y_1 \cdots y_{n+N} |\mathcal{B}(z_1) \cdots \mathcal{B}(z_N) 
|x_1 \cdots x_n \rangle/D
$ where $D:=\prod_{j=1}^N z_j^{j-1-N} (1+t z_j^2)
\prod_{1 \le j < k \le N}(1+t z_j z_k)
(1+t z_j z_k^{-1})$, is invariant under any permutation
$z_i \longleftrightarrow z_j$.

From the fact that the rescaled double row $B$-operator $\mathcal{B}^\prime(z)$
consists of rescaled $L$-operators and $K$-matrix,
which can be essentially obtained from
the original $L$-operators and $K$-matrix
by the changing the spectral parameters $z \longrightarrow t^{-1}z$,
the statement of the invariance
for the matrix element
$ \langle y_1 \cdots y_{n+N} |\mathcal{B}(z_1) \cdots \mathcal{B}(z_N) 
|x_1 \cdots x_n \rangle$
can be converted to
the statement that every matrix element
$\langle \overline{x_1} \cdots \overline{x_n} 
|\mathcal{B}^\prime(z_1) \cdots \mathcal{B}^\prime(z_N) 
|\overline{y_1} \cdots \overline{y_{n+N}} \rangle$
has the property that
$t^N \langle \overline{x_1} \cdots \overline{x_n} 
|\mathcal{B}^\prime(z_1) \cdots \mathcal{B}^\prime(z_N) 
|\overline{y_1} \cdots \overline{y_{n+N}} \rangle/D^\prime$
where $D^\prime=\prod_{j=1}^N z_j^{j-1-N} (1+t^\prime z_j^2)
\prod_{1 \le j < k \le N}(1+t^\prime z_j z_k)
(t^{\prime}+z_j z_k^{-1})$, $t^\prime=t^{-1}$,
is invariant under any permutation $z_i \longleftrightarrow z_j$.
Property 2 follows from this statement since we are considering
a special matrix element
$t^N \langle 1 \cdots M|\mathcal{B}^\prime(z_1)
\cdots \mathcal{B}^\prime(z_N)|\overline{x_1} \cdots
\overline{x_N} \rangle/D^\prime$.

We prove Property 3 as follows.
Since one cannot apply Ivanov's argument as exactly as it is
due to the change of the matrix elements of the $L$-operators,
we modify the argument as follows.
First, we find that it is enough to show
\begin{align}
\langle 1 \cdots M|\mathcal{B}^\prime(z)|
\overline{x} \rangle
=\frac{(1+t^{-1}z^2)}{tz} \frac{z^{\overline{\lambda}+1}-z^{-\overline{\lambda}-1}}{z-z^{-1}},
\ \ \ \overline{\lambda}=\overline{x}-1.
\label{showthislateron}
\end{align}
From
\eqref{showthislateron} and using the
decomposition,
\begin{align}
\langle 1 \cdots M|\mathcal{B}^\prime(z_1)
\cdots \mathcal{B}^\prime(z_N)|\overline{x_1} \cdots
\overline{x_N} \rangle
=\sum_{\overline{y}}
\langle 1 \cdots M|\mathcal{B}^\prime(z_1)|\overline{y} \rangle
\langle \overline{y}|\mathcal{B}^\prime(z_2) \cdots \mathcal{B}^\prime(z_N)|
\overline{x_1} \cdots
\overline{x_N} \rangle,
\end{align}
we have
\begin{align}
\langle 1 \cdots M|\mathcal{B}^\prime(z_1)
\cdots \mathcal{B}^\prime(z_N)|\overline{x_1} \cdots
\overline{x_N} \rangle
=\frac{(1+t^{-1}z_1^2)}{tz_1}
\sum_{\overline{y}}
\frac{z_1^{\overline{y}}-z_1^{-\overline{y}}}{z_1-z_1^{-1}}
\langle \overline{y}|\mathcal{B}^\prime(z_2) \cdots \mathcal{B}^\prime(z_N)|
\overline{x_1} \cdots
\overline{x_N} \rangle,
\end{align}
from which one gets
\begin{align}
\frac{t^N \langle 1 \cdots M|\mathcal{B}^\prime(z_1)
\cdots \mathcal{B}^\prime(z_N)|\overline{x_1} \cdots
\overline{x_N} \rangle}
{
t^N \langle 1 \cdots M|\mathcal{B}^\prime(z_1)
\cdots \mathcal{B}^\prime(z_N)|\overline{x_1} \cdots
\overline{x_N} \rangle|_{z_1 \longleftrightarrow z_1^{-1}}
}
=\frac{(t^\prime z_1+z_1^{-1})}{(z_1+t^\prime z_1^{-1})}.
\end{align}
This ratio and
$\displaystyle
\frac{D^\prime}{D^\prime|_{z_1 \longleftrightarrow z_1^{-1}}}
=\frac{(t^\prime z_1+z_1^{-1})}{(z_1+t^\prime z_1^{-1})}
$ gives the relation
\begin{align}
&t^N \langle 1 \cdots M|\mathcal{B}^\prime(z_1)
\cdots \mathcal{B}^\prime(z_N)|\overline{x_1} \cdots
\overline{x_N} \rangle/D^\prime \nonumber \\
=&t^N \langle 1 \cdots M|\mathcal{B}^\prime(z_1)
\cdots \mathcal{B}^\prime(z_N)|\overline{x_1} \cdots
\overline{x_N} \rangle/D^\prime \Big|_{z_1 \longleftrightarrow z_1^{-1}},
\end{align}
which shows that
$t^N \langle 1 \cdots M|\mathcal{B}^\prime(z_1)
\cdots \mathcal{B}^\prime(z_N)|\overline{x_1} \cdots
\overline{x_N} \rangle/D^\prime$
is invariant under $z_1 \longleftrightarrow z_1^{-1}$.
Together with Property 2, one gets Property 3.

Let us show \eqref{showthislateron}.
We insert the completeness relation to decompose
$
\langle 1 \cdots M|\mathcal{B}^\prime(z)|
\overline{x} \rangle$
as
\begin{align}
&\langle 1 \cdots M|\mathcal{B}^\prime(z)|
\overline{x} \rangle
\nonumber \\
=&
t^{-1}z \sum_{\overline{y}}
\langle 1 \cdots M| \widetilde{B}^\prime(z) | \overline{y} \rangle
\langle \overline{y} |A^\prime(z)|\overline{x} \rangle
+z^{-1}
\langle 1 \cdots M| \widetilde{A}^\prime(z) | 1 \cdots M \rangle
\langle 1 \cdots M|B^\prime(z)|\overline{x} \rangle
\nonumber \\
=&
t^{-1}z \sum_{\overline{y}=1,\dots,\overline{x}-1}
\langle 1 \cdots M| \widetilde{B}^\prime(z) | \overline{y} \rangle
\langle \overline{y} |A^\prime(z)|\overline{x} \rangle
\nonumber \\
&+t^{-1}z
\langle 1 \cdots M| \widetilde{B}^\prime(z) | \overline{x} \rangle
\langle \overline{x} |A^\prime(z)|\overline{x} \rangle
+z^{-1}
\langle 1 \cdots M| \widetilde{A}^\prime(z) | 1 \cdots M \rangle
\langle 1 \cdots M|B^\prime(z)|\overline{x} \rangle.
\label{usetherighthandside}
\end{align}
Inserting the explicit forms of the matrix elements
\begin{align}
\langle 1 \cdots M| \widetilde{B}^\prime(z) | \overline{y} \rangle
&=z^{\overline{y}-1}, \ \ \ \overline{y}=1,\dots,\overline{x}-1, \\
\langle \overline{y} |A^\prime(z)|\overline{x} \rangle
&=\frac{t+1}{t} z^{-\overline{x}+\overline{y}}, \ \ \ \overline{y}=1,\dots,\overline{x}-1, \\
\langle 1 \cdots M| \widetilde{B}^\prime(z) | \overline{x} \rangle
&=z^{\overline{x}-1}, \\
\langle \overline{x} |A^\prime(z)|\overline{x} \rangle
&=\frac{1}{t}, \\
\langle 1 \cdots M| \widetilde{A}^\prime(z) | 1 \cdots M \rangle
&=1, \\
\langle 1 \cdots M|B^\prime(z)|\overline{x} \rangle
&=\frac{1}{t}z^{-\overline{x}+1},
\end{align}
into the right hand side of
\eqref{usetherighthandside}, we have
\begin{align}
\langle 1 \cdots M|\mathcal{B}^\prime(z)|
\overline{x} \rangle
=&t^{-1}z \sum_{\overline{y}=1,\dots,\overline{x}-1}
\frac{t+1}{t} z^{-\overline{x}+\overline{y}} z^{\overline{y}-1}
+t^{-1}z \frac{1}{t}z^{\overline{x}-1}
+z^{-1}\frac{1}{t}z^{-\overline{x}+1} \nonumber \\
=&t^{-2}(t+1) \frac{z^{\overline{x}-1}-z^{1-\overline{x}}}{z-z^{-1}}
+t^{-2}z^{\overline{x}}+t^{-1}z^{-\overline{x}} \nonumber \\
=&\frac{z+tz^{-1}}{t^2(z-z^{-1})}(z^{\overline{x}}-z^{-\overline{x}})
=
\frac{(1+t^{-1}z^2)}{tz} \frac{z^{\overline{\lambda}+1}-z^{-\overline{\lambda}-1}}{z-z^{-1}},
\end{align}
and
\eqref{showthislateron} is shown.

Properties 1, 2 and 3 show that
$t^N \langle 1 \cdots M|\mathcal{B}^\prime(z_1)
\cdots \mathcal{B}^\prime(z_N)|\overline{x_1} \cdots
\overline{x_N} \rangle$ is a polynomial of $t^\prime$
with highest degree $2N^2$,
and is divided by $D^\prime=
\prod_{j=1}^N z_j^{j-1-N} (1+t^\prime z_j^2)
\prod_{1 \le j < k \le N}(1+t^\prime z_j z_k)
(t^{\prime}+z_j z_k^{-1})$.
Hence, Lemma \ref{lemmaone} is proved.

\end{proof}

From Lemma \ref{lemmaone}, one sees that to study the wavefunction,
it is enough to examine a particular value of $t$.
At the point $t=-1$, the six-vertex model
reduces to a five-vertex model since the first $L$-operator becomes
\begin{eqnarray}
L^\prime(z,-1)=
-L(-z,-1)=
\left( 
\begin{array}{cccc}
-1 & 0 & 0 & 0 \\
0 & 1 & -1 & 0 \\
0 & 0 & z^{-1} & 0 \\
0 & 0 & 0 & z^{-1}
\end{array}
\right),
\end{eqnarray}
and second $L$-operator becomes
\begin{eqnarray}
\widetilde{L}^\prime(z,-1)
=\widetilde{L}(-z,-1)
=\left( 
\begin{array}{cccc}
-z & 0 & 0 & 0 \\
0 & z & 1 & 0 \\
0 & 0 & 1 & 0 \\
0 & 0 & 0 & 1
\end{array}
\right).
\end{eqnarray}
The $K$-matrix at $t=-1$ becomes
\begin{eqnarray}
K^\prime(z,-1)=-K(-z,-1)=\left( 
\begin{array}{cc}
tz & 0 \\
0 & z^{-1} \\
\end{array}
\right).
\end{eqnarray}
It is easy to examine at this reduced point $t=-1$,
and we find the following relation.
\begin{lemma}
We have
\begin{align}
&\prod_{j=1}^N z_j^{N+1-j} (1+t^{-1}z_j^2)^{-1}
\prod_{1 \le j<k \le N}(1+t^{-1}z_j z_k)^{-1}
(t^{-1}+z_j z_k^{-1})^{-1} \nonumber \\
\times&t^N
\langle 1 \cdots M|\mathcal{B}^\prime(z_1) \cdots
\mathcal{B}^\prime(z_N)| \overline{x_1} \cdots \overline{x_N}
\rangle \Bigg|_{t=-1}
=sp_{\overline{\lambda}}(\{ z \}_N). \label{forproof}
\end{align}
\end{lemma}
\begin{proof}
From the factorization formula of the determinant
\begin{align}
\mathrm{det}_N(z_j^{N-k+1}-z_j^{-N+k-1})
=(-1)^N \prod_{j=1}^N z_j^{j-1-N}(1-z_j^2)
\prod_{1 \le j <k \le N}(1-z_j z_k)(1-z_j z_k^{-1}),
\end{align}
and the definition of the symplectic Schur functions \eqref{symplecticSchur},
one sees that proving the Lemma is equivalent to showing
the following equality
\begin{align}
\langle 1 \cdots M|\mathcal{B}^\prime(z_1) \cdots
\mathcal{B}^\prime(z_N)| \overline{x_1} \cdots \overline{x_N}
\rangle=(-1)^{N(N-1)/2}
\mathrm{det}_N(z_j^{\lambda_k+N-k+1}-z_j^{-\lambda_k-N+k-1}).
\label{whatwewillshow}
\end{align}

To show this, let us first list the matrix elements
of the single $A$- and $B$-operators.
\begin{lemma}
We have the following explicit expressions for the matrix elements
of the $A$- and $B$-operators at $t=-1$. \\
(1) The matrix element of $A^\prime(z)$ is given by
\begin{align}
\langle \overline{x_1} \cdots \overline{x_{k}}
|A^\prime(z)| \overline{y_1}
\cdots \overline{y_{k}} \rangle
=&(-1)^k \prod_{j=1}^k \delta_{\overline{x_j} \overline{y_j}}.
\label{matrixelementst=-1no1}
\end{align}
(2) The matrix element of $B^\prime(z)$ is given by
\begin{align}
\langle \overline{x_1} \cdots \overline{x_{k-1}}
|B^\prime(z)| \overline{y_1}
\cdots \overline{y_{k}} \rangle
=&(-1)^k (-1)^{j-1} z^{1-\overline{y_j}},
\label{matrixelementst=-1no2}
\end{align}
when the hole configurations $\{ \overline{x} \}$ and $\{ \overline{y} \}$
satisfy \\
$\overline{x_1}=\overline{y_1}, \cdots,
\overline{x_{j-1}}=\overline{y_{j-1}}$,
$\overline{x_{j}}=\overline{y_{j+1}}, \cdots,
\overline{x_{k-1}}=\overline{y_{k}}$ for some $j$,
and 0 otherwise. \\
\\
(3) The matrix element of $\widetilde{A}^\prime(z)$ is given by
\begin{align}
\langle \overline{x_1} \cdots \overline{x_{k}}
|\widetilde{A}^\prime(z)| \overline{y_1}
\cdots \overline{y_{k}} \rangle
=&\prod_{j=1}^k \delta_{\overline{x_j} \overline{y_j}},
\label{matrixelementst=-1no3}
\end{align}
(4) The matrix element of $\widetilde{B}^\prime(z)$ is given by
\begin{align}
\langle \overline{x_1} \cdots \overline{x_{k-1}}
|\widetilde{B}^\prime(z)| \overline{y_1}
\cdots \overline{y_{k}} \rangle
=&(-1)^{j-1} z^{\overline{y_j}-1},
\label{matrixelementst=-1no4}
\end{align}
when the hole configurations $\{ \overline{x} \}$ and $\{ \overline{y} \}$
satisfy \\
$\overline{x_1}=\overline{y_1}, \cdots,
\overline{x_{j-1}}=\overline{y_{j-1}}$,
$\overline{x_{j}}=\overline{y_{j+1}}, \cdots,
\overline{x_{k-1}}=\overline{y_{k}}$ for some $j$,
and 0 otherwise.
\end{lemma}

Using these explicit forms of the matrix elements
\eqref{matrixelementst=-1no1},
\eqref{matrixelementst=-1no2},
\eqref{matrixelementst=-1no3}
\eqref{matrixelementst=-1no4},
and using the decomposition of the double row $B$-operator
\eqref{doublerow2},
one finds that the matrix elements of a single double row $B$-operator
$\mathcal{B}^\prime(z)$ at $t=-1$ is given by
\begin{align}
\langle \overline{x_1} \cdots \overline{x_{k-1}}
|\mathcal{B}^\prime(z)| \overline{y_1}
\cdots \overline{y_{k}} \rangle
=&(-1)^{k+1}(-1)^{j-1} (z^{\overline{y_j}}-z^{-\overline{y_j}}),
\label{matrixelementst=-1}
\end{align}
when the hole configurations $\{ \overline{x} \}$ and $\{ \overline{y} \}$
satisfy \\
$\overline{x_1}=\overline{y_1}, \cdots,
\overline{x_{j-1}}=\overline{y_{j-1}}$,
$\overline{x_{j}}=\overline{y_{j+1}}, \cdots,
\overline{x_{k-1}}=\overline{y_{k}}$ for some $j$,
and 0 otherwise.

Since the matrix elements of a single $B$-operator are essentially the same
with the ones for the original wavefunction at $t=-1$ in \cite{Iv}
except the sign $(-1)^{k+1}$ 
(we also have to translate the hole configurations to
particle configurations), the same argument can be applied,
and one finds the wavefunction at $t=-1$ is the symplectic Schur functions
$\mathrm{det}_N(z_j^{\lambda_k+N-k+1}-z_j^{\lambda_k-N+k-1})$
multiplied by a sign factor
$\prod_{k=1}^N (-1)^{k+1}=(-1)^{N(N-1)/2}$.

\begin{align}
&\langle 1 \cdots M|\mathcal{B}^\prime(z_1) \cdots
\mathcal{B}^\prime(z_N)| \overline{x_1} \cdots \overline{x_N}
\rangle|_{t=-1} \nonumber \\
=&
\prod_{k=1}^N (-1)^{k+1}
\sum_{\sigma \in S_N} (-1)^\sigma
\prod_{j=1}^N 
(
z_j^{\overline{x_{\sigma(j)}}}
-
z_j^{-\overline{x_{\sigma(j)}}}
)
\nonumber \\
=&
\prod_{k=1}^N (-1)^{N(N-1)/2}
\sum_{\sigma \in S_N} (-1)^\sigma
\prod_{j=1}^N 
(
z_j^{\overline{\lambda_{\sigma}(j)}+N-\sigma(j)+1}
-
z_j^{-\overline{\lambda_{\sigma}(j)}-N+\sigma(j)-1}
)
\nonumber \\
=&(-1)^{N(N-1)/2}
\mathrm{det}_N(z_j^{\lambda_k+N-k+1}-z_j^{-\lambda_k-N+k-1})
.
\end{align}
Thus \eqref{whatwewillshow} is shown, hence \eqref{forproof} is proved.
\end{proof}
Finally, from Lemma \ref{lemmaone} and \eqref{forproof},
we have
\begin{align}
&\prod_{j=1}^N z_j^{N+1-j} (1+t^{-1}z_j^2)^{-1}
\prod_{1 \le j<k \le N}(1+t^{-1}z_j z_k)^{-1}
(t^{-1}+z_j z_k^{-1})^{-1} \nonumber \\
\times&t^N
\langle 1 \cdots M|\mathcal{B}^\prime(z_1) \cdots
\mathcal{B}^\prime(z_N)| \overline{x_1} \cdots \overline{x_N}
\rangle \nonumber \\
=&\prod_{j=1}^N z_j^{N+1-j} (1+t^{-1}z_j^2)^{-1}
\prod_{1 \le j<k \le N}(1+t^{-1}z_j z_k)^{-1}
(t^{-1}+z_j z_k^{-1})^{-1} \nonumber \\
\times&t^N
\langle 1 \cdots M|\mathcal{B}^\prime(z_1) \cdots
\mathcal{B}^\prime(z_N)| \overline{x_1} \cdots \overline{x_N}
\rangle \Bigg|_{t=-1} \nonumber \\
=&sp_{\overline{\lambda}}(\{ z \}_N).
\end{align}
which is exactly
\eqref{equivalenttwo}, hence
Theorem \ref{theoremdualandsymplecticschur} is proved.
\end{proof}

\section{A generalization to the factorial symplectic Schur functions}
We have showed Theorem \ref{theoremdualandsymplecticschur}
which gives the relation between the dual wavefunction
and the symplectic Schur functions.
The proof given in Ivanov \cite{Iv} and the last section can be lifted to
give the exact correspondence between the wavefunction,
the dual wavefunction and the
factorial symplectic Schur functions
by introducing inhomogeneous parameters
in the quantum spaces.
We state the correspondence in this section.

First we introduce the $L$-operator of Bump-McNamara-Nakasuji \cite{BMN} which
now has dependence on the quantum space $\mathcal{F}_j$.
The $L$-operator $L_{aj}(z,t,\alpha_j)$
at the $j$-th site
in the quantum space is given by
\begin{eqnarray}
L_{aj}(z,t,\alpha_j)=\left( 
\begin{array}{cccc}
1 & 0 & 0 & 0 \\
0 & t & 1 & 0 \\
0 & (t+1)z^{-1} & z^{-1}+\alpha_j & 0 \\
0 & 0 & 0 & z^{-1}-t \alpha_j
\end{array}
\right). \label{generalizedfelderhofloperator}
\end{eqnarray}
The $L$-operators
now have inhomogeneous parameters
$\alpha_j$, $j=1,\cdots M$ besides the spectral parameters
and the deformation parameter.
For the case of the wavefunction without reflecting boundary,
these newly introduced inhomogeneous parameters
become factorial parameters of the factorial Schur functions in the end.

We also introduce inhomogeneous parameters
in the second $L$-operator.
The second $L$-operator $\widetilde{L}_{aj}(z,t,\alpha_j)$
at the $j$-th site
in the quantum space is given by
\begin{eqnarray}
\widetilde{L}_{aj}(z,t,\alpha_j)=\left( 
\begin{array}{cccc}
z+\alpha_j & 0 & 0 & 0 \\
0 & tz-\alpha_j & 1 & 0 \\
0 & (t+1)z & 1 & 0 \\
0 & 0 & 0 & 1
\end{array}
\right). \label{generalizedsecondloperator}
\end{eqnarray}
One can also generalize the $K$-matrix
to the following one
\begin{eqnarray}
K_{a}(z,t,\alpha_0)=\left( 
\begin{array}{cc}
tz-\alpha_0 & 0 \\
0 & z^{-1}+\alpha_0 \\
\end{array}
\right), \label{generalizedkmatrix}
\end{eqnarray}
where $\alpha_0$ is a free parameter.

Using these inhomogeneous $L$-operators and $K$-matrix,
we as again introduce two types of monodromy matrices
\begin{align}
T_{a}(z,\{ \alpha \})=L_{a M}(z,t,\alpha_M) \cdots L_{a 1}(z,t,\alpha_1)
,
\label{generalizedmonodromy1}
\end{align}
and
\begin{align}
\widetilde{T}_{a}(z,\{ \alpha \})=\widetilde{L}_{a 1}(z,t,\alpha_1) \cdots \widetilde{L}_{a M}(z,t,\alpha_M)
,
\label{generalizedmonodromy2}
\end{align}
which act on $W_a \otimes (\mathcal{F}_1\otimes\dots\otimes 
\mathcal{F}_{M})$,
and denote the matrix elements of the two monodromy matrices as
\begin{align}
A(z,\{ \alpha \})={}_a \langle 0|T_{a}(z, \{ \alpha \})|0 \rangle_{a}, \\
B(z,\{ \alpha \})={}_a \langle 0|T_{a}(z, \{ \alpha \})|1 \rangle_{a},
\end{align}
and
\begin{align}
\widetilde{A}(z,\{ \alpha \})={}_a \langle 1|\widetilde{T}_{a}(z,\{ \alpha \})|1 \rangle_{a}, \\
\widetilde{B}(z,\{ \alpha \})={}_a \langle 0|\widetilde{T}_{a}(z,\{ \alpha \})|1 \rangle_{a}.
\end{align}
Here, $\{ \alpha \}=\{\alpha_1,\dots,\alpha_M \}$
is included in the notation to indicitate that the operators depend
on this set of parameters.
As again,
we introduce the following double row $B$-operator.
\begin{align}
\mathcal{B}(z,\{ \overline{\alpha} \})&=\widetilde{B}(z,\{\alpha \}) {}_a \langle 0| K(z,t,\alpha_0)|0 \rangle_a A(z)
+\widetilde{A}(z,\{ \alpha \}) {}_a \langle 1| K(z,t,\alpha_0)|1 \rangle_a B(z,t,\{ \alpha \}) \\
&=(tz-\alpha_0) \widetilde{B}(z,\{ \alpha \})A(z,\{ \alpha \})
+(z^{-1}+\alpha_0) \widetilde{A}(z,\{ \alpha \}) B(z,\{ \alpha \}).
\label{generalizeddoublerow}
\end{align}
Since the double row $B$-operator uses the generalized $K$-matrix
as a component, the dependence on the inhomogeneous parameters
is lifted to the set of parameters
$\{ \overline{\alpha} \}=\{\alpha_0,\alpha_1,\dots,\alpha_M \}$
where $\alpha_0$ is added to $\{ \alpha \}$.

We now introduce the inhomogeneous wavefunction
$\langle x_1 \cdots x_N|\Psi(z_1,\dots,z_N,\{ \alpha \}) \rangle$
as the overlap between the particle configurations
$\langle x_1 \cdots x_N|$
and the inhomogeneous $N$-particle state
\begin{align}
\Psi(z_1,\dots,z_N,\{ \overline{\alpha} \}) \rangle
=
\mathcal{B}(z_1,\{ \overline{\alpha} \}) \cdots
\mathcal{B}(z_N,\{ \overline{\alpha} \})| \Omega \rangle.
\end{align}
Likewise, the dual wavefunction
$\langle \Phi(z_1,\dots,z_N,\{ \alpha \})|
\overline{x_1} \cdots \overline{x_N} \rangle$
is introduced as the overlap between the hole configurations
$|\overline{x_1} \cdots \overline{x_N} \rangle$
and the inhomogeneous dual $N$-particle state
\begin{align}
\langle \Phi(z_1,\dots,z_N,\{ \overline{\alpha} \})|
=\langle 1 \cdots M|
\mathcal{B}(z_1,\{ \overline{\alpha} \}) \cdots
\mathcal{B}(z_N,\{ \overline{\alpha} \}).
\end{align}

These wavefunctions can be expressed by
the factorial symplectic Schur functions defined below.

\begin{definition}
The factorial symplectic Schur functions is defined
to be the following determinant:
\begin{align}
sp_\lambda(\{ z\}_N|\{ \overline{\alpha} \})
=\frac{G_{\lambda+\delta}(\{z\}_N|\{ \overline{\alpha} \})}
{\mathrm{det}_N(z_j^{N-k+1}-z_j^{-N+k-1})}
, \label{factorialsymplecticSchurfunction}
\end{align}
where $\{ z \}_N=\{z_1,\dots,z_N \}$ is a set of variables
and $\lambda$ denotes a Young diagram
$\lambda=(\lambda_1,\lambda_2,\dots,\lambda_N)$
with weakly decreasing non-negative integers
$\lambda_1 \ge \lambda_2 \ge \cdots \ge \lambda_N \ge 0$,
and $\delta=(N-1,N-2,\dots,0)$.
$G_{\mu}(\{ z \}_N|\{ \overline{\alpha} \})$
is an $N \times N$ determinant
\begin{align}
G_{\mu}(\{ z \}_N|\{ \overline{\alpha} \})
=\mathrm{det}_N
\Bigg(
\prod_{j=0}^\mu (z_k+\alpha_j)
-
\prod_{j=0}^\mu (z_k^{-1}+\alpha_j)
\Bigg).
\end{align}
We remark that one must respect the ordering of the factorial parameters
$\{ \overline{\alpha} \}=\{\alpha_0,\alpha_1,\dots,\alpha_M \}$.
\end{definition}
We have the following correspondence
between the wavefunction of the  free-fermion model with
inhomogeneties and the factorial Schur symplectic functions.
\begin{theorem}
The wavefunction
$\langle x_1 \dots x_N|\Psi(z_1,\dots,z_N,\{ \overline{\alpha} \}) \rangle$
is expressed by the factorial symplectic functions as
\begin{align}
&\bra x_1 \cdots x_N | \Psi(z_1,\dots,z_N,\{ \overline{\alpha} \}) \ket
\nonumber \\
=&\prod_{j=1}^N z_j^{j-1-N}(1+tz_j^2)
\prod_{1 \le j<k \le N}(1+tz_j z_k)(1+tz_j z_k^{-1})
sp_\lambda(\{ z \}_N|\{ \overline{\alpha} \}),
\end{align}
under the relation $\lambda_j=x_{N-j+1}-N+j-1$, $j=1,\dots,N$.
\end{theorem}
This Theorem can be proved by noting that the arguments in Ivanov \cite{Iv}
naturally lift to this inhomogeneous setting.
One first shows that the wavefunction is a polynomial
of $t$ with highest degree $N^2$
whose $t$-dependent part can be factorized as
$\prod_{j=1}^N (1+tz_j^2)
\prod_{1 \le j<k \le N}(1+tz_j z_k)(1+tz_j z_k^{-1})$.
Then one evaluates the wavefunction at $t=-1$,
at which the six-vertex model reduces to a five-vertex model,
and each configuration making non-zero contribution
($2^N \times N!$ configurations in total)
to the wavefunction essentially corrresponds to each term
of the determinant expansion of the numerator of the
factorial symplectic Schur functions \eqref{factorialsymplecticSchurfunction}.

As for the inhomogeneous dual wavefunction,
one can apply the argument in section 4 and get the
following relation with the
factorial symplectic Schur functions.
\begin{theorem} \label{theoremdualandfactorialschur}
The dual wavefunction
$\langle \Phi(z_1,\dots,z_N,\{ \overline{\alpha} \})| \overline{x_1} \cdots \overline{x_N}
\rangle$
can be expressed by the factorial symplectic Schur functions as
\begin{align}
&\langle \Phi(z_1,\dots,z_N,\{ \overline{\alpha} \})| \overline{x_1} \cdots \overline{x_N}
\rangle \nonumber \\
=&t^{N(M-N)}
\prod_{j=1}^N z_j^{j-1-N}(1+tz_j^2)
\prod_{1 \le j<k \le N}(1+tz_j z_k)(1+tz_j z_k^{-1})
sp_{\overline{\lambda}} ( \{ tz \}_N |\{-\overline{\alpha} \}).
\label{dualwavefunctionsandfactorialsymplecticschur}
\end{align}
Here the Young diagram for the factorial symplectic Schur functions
corresponds to the hole configuration under the relation
$\overline{\lambda_j}=\overline{x_{N-j+1}}-N+j-1$, $j=1,\dots,N$,
and the symmetric variables are
$\displaystyle \{ tz \}_N=
\{ tz_1,\dots,tz_N \}$.
Moreover, the signs of the parameters of the factorial symplectic
Schur functions in the right hand side of
\eqref{dualwavefunctionsandfactorialsymplecticschur} are now
inverted simultaneously: $\{-\overline{\alpha} \}
=\{-\alpha_0,-\alpha_1,\dots,-\alpha_M \}$.
\end{theorem}
The correspondence \eqref{dualwavefunctionsandfactorialsymplecticschur}
can be proved by naturally lifting the arguments given in section 4
to this inhomogeneous setting.
At the point when $t=-1$ where the six-vertex model reduces to
the five-vertex model, the introduction of inhomogeneous parameters
is reflected in the $t$-independent part of the correspondence.
The right hand side of the final expression of the correspondence in
\eqref{dualwavefunctionsandfactorialsymplecticschur}
is lifted to the factorial symplectic Schur functions.

\section{Conclusion}
We investigated the free-fermion model under the
reflecting boundary condition, and
showed the precise relation between the dual
wavefunction and the symplectic Schur functions.
The result and the proof is an extension of the
ones in \cite{LMP} for the case without reflecting boundary,
where the statement was transformed into another
equivalent one so that one can use the arguments
given by \cite{BBF} and \cite{Iv}.
The correspondence can be regarded as a type $C$ version of the
dual version of the Tokuyama formula by Ivanov.

In its relation with authomorphic representation theory,
the wavefunction with another boundary $K$-matrices
is introduced by Brubaker-Bump-Chinta-Gunnells \cite{BBCG}.
The conjectures about the correpondence
may be solved by using other ideas such as the
theory of divided difference operators.
As for the wavefunction with reflecting boundary condition,
we remark that there are several works on the XXZ model
with reflecting boundary condition and its degeneration
in \cite{Tsuchiya,RK,CMRV,WZnew,vDE}, for example.

We also generalized the correspondence between the wavefunction,
the dual wavefunction and the symplectic Schur functions
to factorial symplectic Schur functions
by using the first inhomogeneous $L$-operator in \cite{BMN},
the second inhomogeneous $L$-operator and the inhomogeneous
$K$-matrix.
The result is a symplectic version of the result in \cite{BMN},
where the correspondence between the wavefunction without reflecting boundary
and the factorial Schur polynomials is established.
We extended furthermore
to the free-fermion model with two types of
inhomogeneous parameters, and there are
correspondences between the original and the dual
wavefunctions and a generalization
of the factorial Schur polynomials
and factorial symplectic Schur functions \cite{MSW1,MSW2}.
Details will appear elsewhere.
Recently, the Tokuyama-type formula for classical groups
were realized combinatorially using the methods of
non-intersecting lattice paths
by Hamel-King \cite{HK,HKFPSAC}. Factorial characters also appear
in their work.
It is worthwhile studying the relation with integrable models,
which may lead to further studies on the subject.

We finally remark that number theorists
regard the free-fermion six-vertex model
as a special case of the ``metaplectic ice"
(see \cite{Meta} for example).
Recently, they established the relation
between the Yang-Baxter equation for the Perk-Schultz model
and the metaplectic ice \cite{BBB}.
It seems worthwhile to study these models
and find novel combinatorial formulas
by means of modern statistical physical methods
and techniques developed to analyze quantum integrable models.

\section*{Acknowledgements}
This work was partially supported by grant-in-Aid
for Research Activity start-up No. 15H06218
and Scientific Research (C) No. 16K05468.

\appendix
\def\thesection{\Alph{section}}
\def\reference{\relax\refpar}

\section{Matrix elements}
We first list the matrix elements of the $A$- and $B$-operators.

\begin{proposition} \label{propositionmatrixelements}
(1)
The matrix elements of a single $A$-operator $A(z)$ is given by
\begin{align}
\langle \overline{x_1} \cdots \overline{x_N} |A(z)| \overline{y_1}
\cdots \overline{y_{N}} \rangle
=&(t+1)^{|\{ \overline{x_j}, \ j=1,\cdots,N \ | \ 
\overline{x_j} \neq \overline{y_{j}}, \ \overline{x_j}
\neq \overline{y_{j-1}} \}|} \nonumber \\
&\times t^{\sum_{j=0}^{N} \mathrm{Max}(\overline{x_{j+1}}
-\overline{y_j}-1, \ 0)}
z^{\sum_{j=1}^N(\overline{x_j}-\overline{y_{j}})},
\label{matrixelementscoordinates}
\end{align}
for hole configurations $\{ \overline{x} \}$
$(1 \le \overline{x_1} < \cdots < \overline{x_N} \le M)$ and
$\{ \overline{y} \}$
$(1 \le \overline{y_1} < \cdots < \overline{y_{N}} \le M)$
satisfying the interlacing relation
$\overline{x_1} \le \overline{y_1} \le \overline{x_2} \le \cdots \le \overline{x_N} \le \overline{y_{N}}$,
and 0 otherwise.
Here we also set $\overline{y_0}=0$ and $\overline{x_{N+1}}=M+1$.
\\
(2)
The matrix elements of a single $B$-operator $B(z)$ is given by
\begin{align}
\langle \overline{x_1} \cdots \overline{x_N} |B(z)| \overline{y_1}
\cdots \overline{y_{N+1}} \rangle
=&(t+1)^{|\{ \overline{x_j}, \ j=1,\cdots,N \ | \ 
\overline{x_j} \neq \overline{y_{j}}, \ \overline{x_j}
\neq \overline{y_{j+1}} \}|} \nonumber \\
&\times t^{\sum_{j=1}^{N+1} \mathrm{Max}(\overline{x_{j}}
-\overline{y_j}-1, \ 0)}
z^{\sum_{j=1}^N(\overline{x_j}-\overline{y_{j+1}})},
\label{matrixelementscoordinates2}
\end{align}
for hole configurations $\{ \overline{x} \}$
$(1 \le \overline{x_1} < \cdots < \overline{x_N} \le M)$ and
$\{ \overline{y} \}$
$(1 \le \overline{y_1} < \cdots < \overline{y_{N+1}} \le M)$
satisfying the interlacing relation
$\overline{y_1} \le \overline{x_1} \le \overline{y_2} \le \overline{x_2} \le \cdots \le \overline{x_N} \le \overline{y_{N+1}}$,
and 0 otherwise.
Here we also set $\overline{x_{N+1}}=M+1$. \\
(3)
The matrix elements of a single $A$-operator $\widetilde{A}(z)$ is given by
\begin{align}
\langle \overline{x_1} \cdots \overline{x_N} |\widetilde{A}(z)| \overline{y_1}
\cdots \overline{y_{N}} \rangle
=&(t+1)^{|\{ \overline{x_j}, \ j=1,\cdots,N \ | \ 
\overline{x_j} \neq \overline{y_{j}}, \ \overline{x_j}
\neq \overline{y_{j-1}} \}|} \nonumber \\
&\times t^{\sum_{j=1}^{N} \mathrm{Max}(\overline{y_{j}}
-\overline{x_j}-1, \ 0)}
z^{\sum_{j=1}^N(\overline{y_j}-\overline{x_{j}})},
\label{matrixelementscoordinates3}
\end{align}
for hole configurations $\{ \overline{x} \}$
$(1 \le \overline{x_1} < \cdots < \overline{x_N} \le M)$ and
$\{ \overline{y} \}$
$(1 \le \overline{y_1} < \cdots < \overline{y_{N}} \le M)$
satisfying the interlacing relation
$\overline{x_1} \le \overline{y_1} \le \overline{x_2} \le \cdots \le \overline{x_N} \le \overline{y_{N}}$,
and 0 otherwise.
We also set $\overline{y_0}=0$. \\
(4)
The matrix elements of a single $B$-operator $\widetilde{B}(z)$ is given by
\begin{align}
\langle \overline{x_1} \cdots \overline{x_N} |\widetilde{B}(z)| \overline{y_1}
\cdots \overline{y_{N+1}} \rangle
=&(t+1)^{|\{ \overline{x_j}, \ j=1,\cdots,N \ | \ 
\overline{x_j} \neq \overline{y_{j}}, \ \overline{x_j}
\neq \overline{y_{j+1}} \}|} \nonumber \\
&\times t^{\sum_{j=1}^{N+1} \mathrm{Max}(\overline{y_{j}}
-\overline{x_{j-1}}-1, \ 0)}
z^{\sum_{j=1}^{N+1}(\overline{y_j}-\overline{x_{j-1}})-1},
\label{matrixelementscoordinates4}
\end{align}
for hole configurations $\{ \overline{x} \}$
$(1 \le \overline{x_1} < \cdots < \overline{x_N} \le M)$ and
$\{ \overline{y} \}$
$(1 \le \overline{y_1} < \cdots < \overline{y_{N+1}} \le M)$
satisfying the interlacing relation
$\overline{y_1} \le \overline{x_1} \le \overline{y_2} \le \overline{x_2} \le \cdots \le \overline{x_N} \le \overline{y_{N+1}}$,
and 0 otherwise.
Here we also set $\overline{x_{0}}=0$.

\end{proposition}

\begin{proof}
The matrix elements \eqref{matrixelementscoordinates2}
is essentially calcaluted in \cite{LMP} (we just need to reverse the
signs of the powers of $z$ to get the result \eqref{matrixelementscoordinates2}
because of the change of the definition of the $L$-operator \eqref{loperator}).
The other cases \eqref{matrixelementscoordinates},
\eqref{matrixelementscoordinates3} and \eqref{matrixelementscoordinates4}
can be calculated in the same way as in \cite{LMP}.

Let us show \eqref{matrixelementscoordinates4} for example.
Let us first count the powers of the spectral parameter $z$.
If the hole configurations
$\{ \overline{x} \}$ and $\{ \overline{y} \}$ are fixed and satisfies
the interlacing relation
$\overline{y_1} \le \overline{x_1} \le \overline{y_2} \le \overline{x_2} \le \cdots \le \overline{x_N} \le \overline{y_{N+1}}$,
the inner states in the auxiliary space is fixed uniquely,
which is a sequence of $0$'s and $1$'s.
We observe that for each sequence $10 \cdots 01$ of the inner states
in the auxiliary space,
all the matrix elements of the
$L$-operators \eqref{secondloperator} in between contribute to the power $z$,
and gives $z^{\sum_j(\overline{y_{j}}-\overline{x_{j-1}})}$ for
some sum over $j$. Taking all of the $01 \cdots 10$ sequences into account,
we have the factor $z^{\sum_{j=1}^{N+1}
(\overline{y_{j}}-\overline{x_{j-1}})-1}$.
Here, we also take into account the first sequence consisting only of $0$'s 
$0 \cdots 0$,  which contribute to the factor
$z^{\overline{y_1}-\overline{x_0}-1}$.

Let us turn to count the powers of $t+1$ and $t$.
We get a factor $t+1$ for each case when both
$\overline{x_j} \neq \overline{y_j}$ and
$\overline{x_j} \neq \overline{y_{j+1}}$ are satisfied
since the matrix element of the $L$-operator is $[L(z,t)]_{01}^{10}=(t+1)z$
at the $\overline{x_j}$-th site for this case.
One gets $(t+1)^{|\{ \overline{x_j}, \ j=1,\cdots,N \ | \ 
\overline{x_j} \neq \overline{y_{j}}, \ \overline{x_j}
\neq \overline{y_{j+1}} \}|}$ in total.

Next, we count the powers of $t$.
If $\overline{y_j}<\overline{x_j}$ is satisfied,
the matrix elements of the $L$-operators are all
$[L(z,t)]_{01}^{01}=tz$ from the $(\overline{x_{j-1}}+1)$-th site to 
the $(\overline{y_j}-1)$-th site. On the other hand, $[L(z,t)]_{01}^{01}$
does not appear if $\overline{x_{j-1}}=\overline{y_j}$, and there is no
contribution to the power of $t$ for this case.
The contributions from $t$ is given by
$t^{\sum_{j=1}^{N+1} \mathrm{Max}(\overline{y_{j}}
-\overline{x_{j-1}}-1, \ 0)}$.

Having calculated all factors, one finds the matrix elements
are given by
\eqref{matrixelementscoordinates4}

\end{proof}

{\bf Example of $\langle \overline{x_1} \cdots \overline{x_N} |A(z)| \overline{y_1}
\cdots \overline{y_{N}} \rangle$} \\
Let $M=15$, $N=4$, $\overline{x}=(3,5,8,11)$ and
$\overline{y}=(3,6,11,13)$.
We also set $\overline{x_5}=15+1=16$ and $\overline{y_0}=0$.
From
$\mathrm{Max}(\overline{x_1}-\overline{y_0}-1,0)=\mathrm{Max}(3-0-1,0)=2$, $\mathrm{Max}(\overline{x_2}-\overline{y_1}-1,0)=\mathrm{Max}(5-3-1,0)=1$,
$\mathrm{Max}(\overline{x_3}-\overline{y_2}-1,0)=\mathrm{Max}(8-6-1,0)=1$, $\mathrm{Max}(\overline{x_4}-\overline{y_3}-1,0)=\mathrm{Max}(11-11-1,0)=0$, $\mathrm{Max}(\overline{x_5}-\overline{y_4}-1,0)=\mathrm{Max}(16-13-1,0)=2$, 
we have the factor $t^{2+1+1+0+2}=t^6$.
The relations $\overline{y_0} \neq \overline{x_1} = \overline{y_1}$,
$\overline{y_1} \neq \overline{x_2} \neq \overline{y_2}$,
$\overline{y_2} \neq \overline{x_3} \neq \overline{y_3}$,
$\overline{y_3} = \overline{x_4} \neq \overline{y_4}$
give the factor $(t+1)^2$, and we also have the factor $z^{-6}$
from $(\overline{x_1}-\overline{y_1})+
(\overline{x_2}-\overline{y_2})+
(\overline{x_3}-\overline{y_3})+
(\overline{x_4}-\overline{y_4})
=(3-3)+(5-6)+(8-11)+(11-13)=0-1-3-2=-6$.
In total, the right hand side of \eqref{matrixelementscoordinates}
is calculated as $(t+1)^2 t^6 z^{-6}$.
One can check that this matches the left hand side
of \eqref{matrixelementscoordinates}, i.e.,
the matrix elements of the corresponding $A$-operator
by explicit calculation. \\

{\bf Example of $\langle \overline{x_1} \cdots \overline{x_N} |B(z)| \overline{y_1}
\cdots \overline{y_{N+1}} \rangle$} \\
Let $M=10$, $N=2$, $\overline{x}=(3,6)$ and
$\overline{y}=(1,6,8)$.
We also set $\overline{x_3}=10+1=11$.
From $\mathrm{Max}(\overline{x_1}-\overline{y_1}-1,0)=\mathrm{Max}(3-1-1,0)=1$, $\mathrm{Max}(\overline{x_2}-\overline{y_2}-1,0)=\mathrm{Max}(6-6-1,0)=0$, $\mathrm{Max}(\overline{x_3}-\overline{y_3}-1,0)=\mathrm{Max}(11-8-1,0)=2$, we have the factor $t^{1+0+2}=t^3$.
The relations $\overline{y_1} \neq \overline{x_1} \neq \overline{y_2}$,
$\overline{y_2} = \overline{x_2} \neq \overline{y_3}$
give the factor $(t+1)^1=t+1$, and we also have the factor $z^{-5}$
from $(\overline{x_1}-\overline{y_2})+
(\overline{x_2}-\overline{y_3})=(3-6)+(6-8)=-3-2=-5$.
In total, the right hand side of \eqref{matrixelementscoordinates2}
is calculated as $(t+1) t^3 z^{-5}$.
One can check that this matches the left hand side
of \eqref{matrixelementscoordinates2}, i.e.,
the matrix elements of the corresponding $B$-operator
by explicit calculation. \\

{\bf Example of $\langle \overline{x_1} \cdots \overline{x_N} |\widetilde{A}(z)| \overline{y_1}
\cdots \overline{y_{N}} \rangle$} \\
Let $M=15$, $N=4$, $\overline{x}=(2,5,10,13)$ and
$\overline{y}=(2,8,10,15)$.
We also set $\overline{y_0}=0$.
From
$\mathrm{Max}(\overline{y_1}-\overline{x_1}-1,0)=\mathrm{Max}(2-2-1,0)=0$, $\mathrm{Max}(\overline{y_2}-\overline{x_2}-1,0)=\mathrm{Max}(8-5-1,0)=2$,
$\mathrm{Max}(\overline{y_3}-\overline{x_3}-1,0)=\mathrm{Max}(10-10-1,0)=0$, $\mathrm{Max}(\overline{y_4}-\overline{x_4}-1,0)=\mathrm{Max}(15-13-1,0)=1$, 
we have the factor $t^{0+2+0+1}=t^3$.
The relations $\overline{y_0} \neq \overline{x_1} = \overline{y_1}$,
$\overline{y_1} \neq \overline{x_2} \neq \overline{y_2}$,
$\overline{y_2} \neq \overline{x_3} = \overline{y_3}$,
$\overline{y_3} \neq \overline{x_4} \neq \overline{y_4}$
give the factor $(t+1)^2$, and we also have the factor $z^{5}$
from $(\overline{y_1}-\overline{x_1})+
(\overline{y_2}-\overline{x_2})+
(\overline{y_3}-\overline{x_3})+
(\overline{y_4}-\overline{x_4})
=(2-2)+(8-5)+(10-10)+(15-13)=0+3+0+2=5$.
In total, the right hand side of \eqref{matrixelementscoordinates3}
is calculated as $(t+1)^2 t^3 z^{5}$.
One can check that this matches the left hand side
of \eqref{matrixelementscoordinates3}, i.e.,
the matrix elements of the corresponding $A$-operator
by explicit calculation. \\

\begin{figure}[h]
\includegraphics[width=15cm]{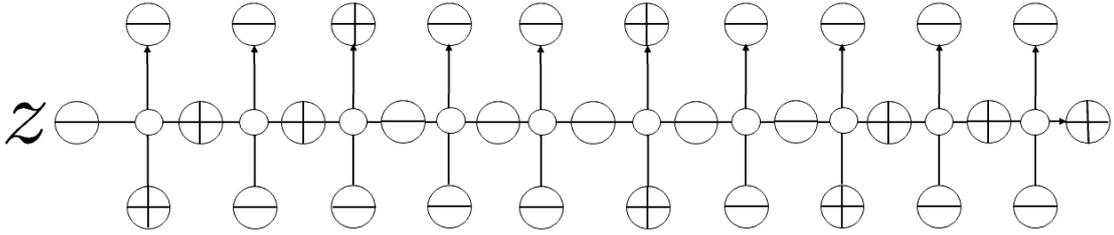}
\caption{The matrix element $\langle \overline{x_1} \cdots \overline{x_N} 
|\widetilde{B}(z)| \overline{y_1}
\cdots \overline{y_{N+1}} \rangle$ for $M=10$, $N=2$,
$\overline{x}=(5,8)$ and
$\overline{y}=(3,5,10)$. One sees that the inner state
is uniquely fixed, and the matrix element is
calculated by multiplying the matrix elements of the $L$-operators
$1 \times tz \times (t+1)z \times 1 \times 1 \times 1 \times 1
\times 1 \times tz \times tz=(t+1) t^3 z^{4}$.}
\label{picturematrixelement}
\end{figure}

{\bf Example of $\langle \overline{x_1} \cdots \overline{x_N} |\widetilde{B}(z)| \overline{y_1}
\cdots \overline{y_{N+1}} \rangle$} \\
Let $M=10$, $N=2$, $\overline{x}=(5,8)$ and
$\overline{y}=(3,5,10)$.
We also set $\overline{x_0}=0$.
From $\mathrm{Max}(\overline{y_1}-\overline{x_0}-1,0)=\mathrm{Max}(3-0-1,0)=2$, $\mathrm{Max}(\overline{y_2}-\overline{x_1}-1,0)=\mathrm{Max}(5-5-1,0)=0$, $\mathrm{Max}(\overline{y_3}-\overline{x_2}-1,0)=\mathrm{Max}(10-8-1,0)=1$, we have the factor $t^{2+0+1}=t^3$.
The relations $\overline{y_1} \neq \overline{x_1} = \overline{y_2}$,
$\overline{y_2} \neq \overline{x_2} \neq \overline{y_3}$
give the factor $(t+1)^1=t+1$, and we also have the factor $z^{4}$
from $(\overline{y_1}-\overline{x_0})+
(\overline{y_2}-\overline{x_1})+(\overline{y_3}-\overline{x_2})-1
=(3-0)+(5-5)+(10-8)-1=3+0+2-1=4$.
In total, the right hand side of \eqref{matrixelementscoordinates4}
is calculated as $(t+1) t^3 z^{4}$.
One can check that this matches the left hand side
of \eqref{matrixelementscoordinates4}, i.e.,
the matrix elements of the corresponding $B$-operator
by explicit calculation
(see Figure \ref{picturematrixelement} for a graphical description
of the corresponding matrix element). \\

Combining the matrix elements of the single $A$- and $B$-operators
\eqref{matrixelementscoordinates},
\eqref{matrixelementscoordinates2},
\eqref{matrixelementscoordinates3} and
\eqref{matrixelementscoordinates4},
one can calculate the matrix elements of the
double row $B$-operators.
\begin{proposition}
The matrix elements of the double row $B$-operator
is given by
\begin{align}
\langle \overline{x^N_N} \cdots \overline{x^N_1} |
\mathcal{B}(z)| \overline{x^{N+1}_1}
\cdots \overline{x^{N+1}_{N+1}} \rangle
=\alpha(z,\{ \overline{x^N} \},\{ \overline{x^{N+1}} \})
+\beta(z,\{ \overline{x^N} \},\{ \overline{x^{N+1}} \}),
\label{doublerowmatrixelements}
\end{align}
for hole configurations $\{ \overline{x^N} \}$
$(1 \le \overline{x^N_1} < \cdots < \overline{x^N_N} \le M)$ and
$\{ \overline{x^{N+1}} \}$
$(1 \le \overline{x^{N+1}_1} < \cdots < \overline{x^{N+1}_{N+1}} \le M)$
satisfying the interlacing relation
$\overline{x^{N+1}_1} \le \overline{x^N_1} \le \overline{x^{N+1}_2} \le \overline{x^N_2} \le \cdots \le \overline{x^N_N} \le \overline{x^{N+1}_{N+1}}$,
and 0 otherwise. \\
$\alpha(z,\{\overline{x^N} \},\{ \overline{x^{N+1}} \})$
is given by
\begin{align}
&\alpha(z,\{\overline{x^N} \},\{ \overline{x^{N+1}} \})
\nonumber \\
=&\sum_{\{ \overline{y^N} \}}
(t+1)^{
|\{ \overline{x^N_j}, \ j=1,\cdots,N \ | \ 
\overline{x^N_j} \neq \overline{y^N_{j}}, \ \overline{x^N_j}
\neq \overline{y^N_{j-1}} \}|
+
|\{ \overline{y^N_j}, \ j=1,\cdots,N \ | \ 
\overline{y^N_j} \neq \overline{x^{N+1}_{j}}, \ \overline{y^N_j}
\neq \overline{x^{N+1}_{j+1}} \}|
} \nonumber \\
&\times t^{
\sum_{j=1}^{N} \mathrm{Max}(\overline{y^N_{j}}
-\overline{x^N_j}-1, \ 0)
+
\sum_{j=1}^{N+1} \mathrm{Max}(\overline{y^N_{j}}
-\overline{x^{N+1}_j}-1, \ 0)
} \nonumber \\
&\times z^{
\sum_{j=1}^N(\overline{y^N_j}-\overline{x^N_{j}})
+\sum_{j=1}^N(\overline{y^N_j}-\overline{x^{N+1}_{j+1}})-1
}, \label{matrixelementspartone}
\end{align}
where we have set $\overline{y^N_{N+1}}=M+1$, $\overline{y^N_0}=0$,
and take the sum over $\{ \overline{y^N} \}
=\{ \overline{y^N_1}, \dots,\overline{y^N_N} \}$
such that $\mathrm{Max}(\overline{x^N_j},\overline{x^{N+1}_j})
\le \overline{y^N_{j}} \le \mathrm{Min}(\overline{x^N_{j+1}},\overline{x^{N+1}_{j+1}})$ is satisfied for each $j=1,\dots,N$.

$\beta(z,\{\overline{x^N} \},\{ \overline{x^{N+1}} \})$
is given by
\begin{align}
&\beta(z,\{\overline{x^N} \},\{ \overline{x^{N+1}} \})
\nonumber \\
=&\sum_{\{ \overline{y^N} \}}
(t+1)^{
|\{ \overline{x^N_j}, \ j=1,\cdots,N \ | \ 
\overline{x^N_j} \neq \overline{y^N_{j}}, \ \overline{x^N_j}
\neq \overline{y^N_{j+1}} \}|
+
|\{ \overline{y^N_j}, \ j=1,\cdots,N+1 \ | \ 
\overline{y^N_j} \neq \overline{x^{N+1}_{j}}, \ \overline{y^N_j}
\neq \overline{x^{N+1}_{j-1}} \}|
} \nonumber \\
&\times t^{
\sum_{j=1}^{N+1} \mathrm{Max}(\overline{y^N_{j}}
-\overline{x^N_{j-1}}-1, \ 0)
+
\sum_{j=0}^{N+1} \mathrm{Max}(\overline{y^N_{j+1}}
-\overline{x^{N+1}_j}-1, \ 0)+1
} \nonumber \\
&\times z^{
\sum_{j=1}^{N+1}(\overline{y^N_j}-\overline{x^N_{j-1}})
+\sum_{j=1}^{N+1}(\overline{y^N_j}-\overline{x^{N+1}_{j}})
}, \label{matrixelementsparttwo}
\end{align}
where we have set $\overline{y^N_{N+2}}=M+1$,
and take the sum over $\{ \overline{y^N} \}
=\{ \overline{y^N_1}, \dots,\overline{y^N_N} \}$
such that $\mathrm{Max}(\overline{x^N_{j-1}},\overline{x^{N+1}_{j-1}})
\le \overline{y^N_{j}} \le \mathrm{Min}(\overline{x^N_{j}},\overline{x^{N+1}_{j}})$ is satisfied for each $j=1,\dots,N$.
Here we also set $\overline{x^N_0}=\overline{x^{N+1}_0}=0$. \\
\end{proposition}

$\alpha(z,\{ \overline{x^N} \},\{ \overline{x^{N+1}} \})$
and $\beta(z,\{ \overline{x^N} \},\{ \overline{x^{N+1}} \})$
are the explicit forms
of the matrix elements
$\langle \overline{x^N_N} \cdots \overline{x^N_1} |
z^{-1} \widetilde{A}(z) B(z)| \overline{x^{N+1}_1}
\cdots \overline{x^{N+1}_{N+1}} \rangle$
and
$\langle \overline{x^N_N} \cdots \overline{x^N_1} |
tz \widetilde{B}(z) A(z)| \overline{x^{N+1}_1}
\cdots \overline{x^{N+1}_{N+1}} \rangle$
which are calculated by inserting the completeness relation
between $\widetilde{A}(z)$ and $B(z)$
and using
\eqref{matrixelementscoordinates2}
\eqref{matrixelementscoordinates3},
and by inserting the completeness relation
between $\widetilde{B}(z)$ and $A(z)$
and using
\eqref{matrixelementscoordinates}
\eqref{matrixelementscoordinates4}, respectively.
The sum of those two matrix elements
becomes the explicit form of the matrix elements
$
\langle \overline{x^N_N} \cdots \overline{x^N_1} |
\mathcal{B}(z)| \overline{x^{N+1}_1}
\cdots \overline{x^{N+1}_{N+1}} \rangle$
due to the decomposition
of the double row $B$-operator \eqref{doublerow}.

Inserting the completeness relation between
each double row $B$-operators and
using \eqref{doublerowmatrixelements} repeatedly,
one gets the explicit form of the dual wavefunction
$\langle \Phi(z_1,\dots,z_N)|\overline{x^N_1} \cdots \overline{x^N_N} \ket$.
\begin{proposition}
The matrix elements of the dual wavefunction
$\langle \Phi(z_1,\dots,z_N)|\overline{x^N_1} \cdots \overline{x^N_N} \ket$
is given by
\begin{align}
\langle \Phi(z_1,\dots,z_N)|\overline{x^N_1} \cdots \overline{x^N_N} \ket
=\sum_{\{ \overline{x^1} \},\dots,\{ \overline{x^{N-1}} \}}
\prod_{k=1}^N
(\alpha(z_k,\{ \overline{x^{k-1}} \},\{ \overline{x^{k}} \})
+\beta(z_k,\{ \overline{x^{k-1}} \},\{ \overline{x^{k}} \})),
\end{align}
where the sum is over all sequences of
hole configurations $\{ \overline{x^k} \}$, $k=1,\dots,N-1$
$(1 \le \overline{x^k_1} < \cdots < \overline{x^k_k} \le M)$
satisfying the interlacing relations
$\overline{x^{k+1}_1} \le \overline{x^{k}_1} \le \overline{x^{k+1}_2} \le \overline{x^{k}_2} \le \cdots \le \overline{x^{k}_{k}} \le \overline{x^{k+1}_{k+1}}$
for all $k=1,\dots,N-1$.
$\alpha(z,\{ \overline{x^{k-1}} \},\{ \overline{x^{k}} \})$
and
$\beta(z,\{ \overline{x^{k-1}} \},\{ \overline{x^{k}} \}))$
are given by
\eqref{matrixelementspartone} and \eqref{matrixelementsparttwo}.
\end{proposition}

Combining the above result with
\eqref{dualwavefunctionsandsymplecticschur} in
Theorem \ref{theoremdualandsymplecticschur},
one gets the following.

\begin{corollary}
The following identity holds
\begin{align}
&t^{N(M-N)}
\prod_{j=1}^N z_j^{j-1-N}(1+tz_j^2)
\prod_{1 \le j<k \le N}(1+tz_j z_k)(1+tz_j z_k^{-1})
sp_{\overline{\lambda}} ( \{ tz \}_N ) \nonumber \\
=&\sum_{\{ \overline{x^1} \},\dots,\{ \overline{x^{N-1}} \}}
\prod_{k=1}^N
(\alpha(z_k,\{ \overline{x^{k-1}} \},\{ \overline{x^{k}} \})
+\beta(z_k,\{ \overline{x^{k-1}} \},\{ \overline{x^{k}} \})),
\end{align}
where $\alpha(z,\{ \overline{x^{k-1}} \},\{ \overline{x^{k}} \})$
and
$\beta(z,\{ \overline{x^{k-1}} \},\{ \overline{x^{k}} \}))$
are given by
\eqref{matrixelementspartone} and \eqref{matrixelementsparttwo},
and
the sum is over all sequences of
hole configurations $\{ \overline{x^k} \}$, $k=1,\dots,N-1$
$(1 \le \overline{x^k_1} < \cdots < \overline{x^k_k} \le M)$
satisfying the interlacing relations
$\overline{x^{k+1}_1} \le \overline{x^{k}_1} \le \overline{x^{k+1}_2} \le \overline{x^{k}_2} \le \cdots \le \overline{x^{k}_{k}} \le \overline{x^{k+1}_{k+1}}$
for all $k=1,\dots,N-1$.
$\overline{x^N_j}$, $j=1,\dots,N$
are uniquely fixed by $\overline{\lambda_j}$, $j=1,\dots,N$
under the relation $\overline{x^N_j}=\overline{\lambda_{N+1-j}}+j$,
$j=1,\dots,N$.
\end{corollary}

\end{document}